\documentclass[11pt]{article}
\usepackage{amsmath}
\usepackage{amssymb}
\usepackage{graphicx, color}

\newtheorem{theorem}{Theorem}[section]
 
 \newtheorem{lemma}{Lemma}[section]
 \newtheorem{proposition}{Proposition}[section]
 
 \newtheorem{remark}{Remark}[section]
 \numberwithin{equation}{section}
 
\def\ds{\displaystyle}
\def\forall{\hbox{for all}~}
\def\TV{\hbox{Tot.Var.}}
\def\L{{\bf L}}

\def\ve{\varepsilon}

\def\D{{\cal D}}
\def\S{{\cal S}}
\def\T{{\cal T}}

\def\R{I\!\!R}

\def\vp{\varphi}

\def\vsk{\vskip 4em}
\def\v{\vskip 1em}
\def\O{{\cal O}}

\def\ov{\overline}
\def\Tilde{\widetilde}
\def\Hat{\widehat}
\def\bega{\begin{array}}
\def\enda{\end{array}}
\def\begi{\begin{itemize}}
\def\endi{\end{itemize}}

\def\bel{\begin{equation}\label}
\def\beq{\begin{equation}}
\def\eeq{\end{equation}}
\def\sqr#1#2{\vbox{\hrule height .#2pt
\hbox{\vrule width .#2pt height #1pt \kern #1pt
\vrule width .#2pt}\hrule height .#2pt }}
\def\square{\sqr74}
\def\endproof{\hphantom{MM}\hfill\llap{$\square$}\goodbreak}

\begin{document}
\title{\bf On finite time BV blow-up for the
p-system}
\vsk

\author{Alberto Bressan$^{(*)}$, Geng Chen$^{(**)}$, and Qingtian Zhang$^{(***)}$\\    \\
(*) Department of Mathematics, Penn State University,
\\
(**) Department of Mathematics, University of Kansas, Lawrence,\\
(***) Department  of Mathematics,
University of California, Davis.
\\
\,\\
e-mails:~ bressan@math.psu.edu,~gengchen@ku.edu,~qzhang@math.ucdavis.edu}

\maketitle

\begin{abstract} The paper studies the possible blowup 
of the total variation
for entropy weak solutions of the  p-system, modeling isentropic gas dynamics. It is assumed that the density remains uniformly positive, 
while the
initial data can have arbitrarily large total variation
(measured in terms of Riemann invariants).
Two main results are proved.  (I) If the total variation
 blows up in finite time,
then the solution must contain an infinite number of large shocks
in a neighborhood of some point in the $t$-$x$ plane.
(II) Piecewise smooth approximate solutions can be constructed whose
total variation blows up in finite time. For these solutions 
the strength of waves
emerging from each interaction is exact, while rarefaction waves
satisfy  the natural decay estimates stemming from the assumption of 
genuine nonlinearity.  
\end{abstract}

\section{Introduction}
\label{sec:0}
\setcounter{equation}{0}

For hyperbolic systems of conservation laws in one space dimension,
a major remaining open problem is whether, for large BV initial data,
the total variation of entropy-weak solutions remains uniformly
bounded or can blow up in finite time.

In the literature, BV bounds have been established
by two main approaches:

{\bf (I)} Estimating the strength of new waves generated at each interaction, regardless of the order in which different wave-fronts cross each other.  For small initial data, this technique was introduced by 
Glimm \cite{G}. Under additional hypothesis, it
can be applied also to solutions with large data. See for example 
\cite{AS, MS, N, NS, Y}.
 
{\bf (II)} Relying on the decay of rarefaction waves, due to genuine nonlinearity, to provide additional cancellations. This approach 
first appeared in \cite{GL} and was then extended in \cite{BCM, C}. 

On the other hand, some particular  $3\times 3$ hyperbolic systems have
been constructed in \cite{BJ, J}, admitting solutions 
whose total variation blows up in finite time.  
One should remark, however, that these
systems do not come from physical models and do not admit a strictly convex entropy.

In this paper we study the possible blowup for solutions to
the p-system
\bel{1}\left\{\begin{array}{rl} v_t - u_x&=~0\,,\cr
u_t + p(v)_x&=~0\,,\end{array}\right.\eeq
modeling  isentropic gas dynamics
in Lagrangian variables.
Here $u$ is the velocity, $\rho$ is the density,  
$v=\rho^{-1}$ is specific volume,
while $p=p(v)$ is the pressure.

In \cite{BCZ, BCZZ}, for a general class of pressure 
functions $p(\cdot)$, the authors  constructed 
piecewise constant  approximate solutions
whose total variation grows without bound.    For these
front tracking approximations, the strength of wave fronts
emerging at each interaction is the same as in the exact solution,
while the only source of error is in the wave speeds.
These examples confirm the  analysis in \cite{CJ}, and
show that uniform BV bounds cannot be established relying only on 
an accurate estimate of wave strengths across interations.

The next question, which we investigate in the present paper, 
is whether BV bounds for the p-system can be established by
taking into account also the decay of rarefaction waves, 
stemming from the assumption  $p''(v)>0$. 
We recall that
Oleinik-type estimates on the decay of positive waves
for genuinely nonlinear $n\times n$ hyperbolic systems
were proved in \cite{Bbook, BC2, BY}.   
The analysis in the last section of \cite{BCZ} shows that,
if this decay of rarefaction waves were taken into account,
then the interaction patterns considered in \cite{BCZ, BCZZ}
would no longer yield a large amplification of total variation.
It is thus natural to ask:

{\bf (Q)} {\it 
Consider a piecewise smooth  approximate solution of (\ref{1}) 
with large BV initial data. Assume that
\begi
\item at each 
interaction, the strengths of outgoing waves are the same 
as in an exact solution,
\item rarefaction waves satisfy decay estimates as in \cite{Bbook, BC2, BY}, due to genuine nonlinearity,
\item the density remains uniformly positive.
\endi
Can the total variation still blow up in finite time?}
\v
An example will be constructed, showing that finite time 
BV blowup for such approximate solutions is indeed possible.   

Although our solutions are not exact, because some errors 
occur in the wave speeds, they possess all the qualitative properties
known for exact solutions.  
The present analysis 
thus provides some indication that finite time blowup 
of the total variation might 
be possible, for the p-system. 

Our second main result yields
 a necessary condition for blowup. 
 Namely, we prove that if the total variation
blows up in finite time, 
then the solution must contain an infinite number of large shocks,
in a neighborhood of some point in the $t$-$x$ plane.

This result should be compared with earlier literature, proving 
BV stability for various classes of initial data
\bel{idata}
u(0,x)~=~\bar u(x),\qquad\qquad v(0,x)~=~\bar v(x)\,.\eeq
\begi
\item If 
$\bar u,\bar v$ have sufficiently small total oscillation,
in the sense that all the initial values 
$(\bar u(x),\bar v(x))$
are contained in a disc with sufficiently small radius,
in the $v$-$u$ plane, then  the solution of (\ref{1}) 
exists globally in time. Moreover, bounds 
on the
total variation can be provided, uniformly in time \cite{GL}.
\item If $\bar u,\bar v$ are a sufficiently small 
BV perturbation of some (possibly large) Riemann data,
then again the solution exists globally in time and its
total variation remains uniformly bounded \cite{BC1, Lw}.
\endi
Building upon these ideas, our present analysis shows
that, for a solution containing only finitely many large shocks,
the total variation remains bounded. 
Indeed,
the blow-up of the BV norm in finite time 
requires the presence of infinitely many large shocks
in a bounded region of the $t$-$x$ plane.

The remainder of the paper is organized as follows. To keep the exposition 
self-contained, in Section 2  we review
some well known results  on the interaction of elementary waves for the p-system. Section~3 develops some estimates related to
the decay of rarefaction waves, valid also for solutions
with large oscillation.
In Section~4 
we construct a piecewise constant approximate solution with a periodic interaction pattern, and where all rarefaction waves
decay at the rate $1/t$.   
By a suitable modification of this basic pattern, 
in Section 5, we construct a piecewise smooth approximate solution 
whose BV norm blows up in finite time.

Section~6 contains the statement of our main theorem, 
providing a necessary condition for finite time blowup.  
Details of the proof are then worked out in Sections~7 and 8.

\section{Elementary wave interactions for the p-system}
\label{sec:1}
\setcounter{equation}{0}
Throughout this paper
 we consider the p-system (\ref{1})  with $\gamma$-law pressure 
\beq\label{2}
p(v)~=~ A v^{-\gamma}~= ~A\rho^\gamma ,
\eeq
for some constants $\gamma>1$ and $A>0$.
For this system one can define 
the Riemann invariants $w_1$ and $w_2$ by setting
\beq\label{RI}
w_1~\doteq~u-h \qquad w_2~\doteq~u+h \,,
\eeq
where
\beq\label{xi_eq}
 h~\doteq~B\,v^{(1-\gamma)/2}~=~B\,\rho^{(\gamma-1)/2}, \quad\qquad  B\doteq \frac2{\gamma-1}\sqrt{A\gamma}\,.
\eeq
For future use, we record the identities 
\bel{i2} 
v~=~\left({h\over B}\right)^{2/(1-\gamma)}~=~\left({(\gamma-1)(w_2-w_1)\over 4\sqrt{A\gamma}}\right)^{2/(1-\gamma)}\,.\eeq

For any smooth solution, these Riemann invariants
remain constant along forward and backward characteristics, respectively.
Namely
\bel{w12}
w_{1,t}-c\,w_{1,x}~=~0\,,\qquad \qquad w_{2,t}+c\,w_{2,x}~=~0\,,
\eeq
where the (Lagrangian) wave speed $c$ is 
\bel{speed}
c~=~\sqrt{-{p}'(v)}~=~\sqrt{A/\gamma}\cdot v^{-(\gamma+1)/2}.
\eeq
In the following, it will be convenient to express the
wave speed in terms of the Riemann coordinates $w_1,w_2$ in 
(\ref{RI}). Introducing the function
\bel{cdef} c(s) ~=~ \sqrt{A/\gamma}\cdot \left({(\gamma-1)\,s\over 4\sqrt{A\gamma}}\right)^{(\gamma+1)/(\gamma-1)}\,,
\eeq
by (\ref{i2}) and (\ref{speed}) the wave speed can be written as
$$\sqrt{-p'(v)}~=~c(w_2-w_1).$$
\v
{\bf Example.} 
In the special case where $p=\rho^3/3$, 
one has the simple relation $h=\rho$.  
By (\ref{RI}) and (\ref{speed}),
the Riemann invariants and the wave speed are given by 
$$w_1~=~u-\rho,\qquad\quad w_2~=~u+\rho,\qquad
\qquad c~=~\rho^2.$$

\subsection{Elementary waves.}
\label{sub_2.1}
A solution to the p-system contains three types of 
 waves: rarefactions, compressions, and shock waves. 
 In terms of the variable $h$ at (\ref{xi_eq}), left and right states will be denoted by 
\bel{uhpm}
(u_-,\, h_-)\qquad\hbox{and} \qquad (u_+, h_+),\eeq 
respectively.

Recalling (\ref{RI}), the signed wave strength will 
always be measured in terms of Riemann invariants:
\bel{strength}
w_{1,+} - w_{1-}~~ \hbox{for a  $1$-wave,}\qquad\quad
w_{2,+} - w_{2,-}~~ \hbox{for a $2$-wave}.
\eeq

We now recall the construction of basic wave curves.  See 
\cite{Bbook, Sm} for details.

\paragraph{1.} 
The rarefaction and compression waves satisfy the following equations:
\begin{itemize}
\item For a 1-wave (backward moving front),
$$u_+-u_-~= ~h_--h_+$$
where $h_+>h_-$ for a 1-compression wave and $h_+<h_-$ for 1-rarefaction wave.

\item For a 2-wave (forward moving front),
$$u_+-u_-~=~ h_+-h_-$$
where $h_+<h_-$ for a 2-compression wave and $h_+>h_-$ for a 2-rarefaction wave.
\end{itemize}

\paragraph{2.} 
A shock wave with left state $(u_-,\rho_- )$ and right state $(u_+,\rho_+)$,
traveling with speed
$\lambda$,
satisfies the Rankine-Hugoniot equations
$$\left\{ \begin{array}{rl} \ds \lambda\left({1\over \rho_+}- {1\over \rho_-}
\right)&=~u_--u_+\,,\cr
&\cr
\lambda(u_+-u_-) & = ~p(v_+) - p(v_-) ~=\ds
~A\rho_+^\gamma - A\rho_-^\gamma\,.\end{array}
\right.$$
The Lax admissibility condition here yields $u_+<u_-\,$ for both 1-waves and 2-waves.
Hence
\bel{RH2}u_+-u_- ~=~- \sqrt{ \left({1\over\rho_+}- {1\over \rho_-}\right)\left(
A\rho_-^\gamma-A\rho_+^\gamma\right)}\,,\eeq
and
\bel{RH5}\lambda ~=~\pm \sqrt A
{\ds\sqrt{\ds {\rho^\gamma_-}-{\rho^\gamma_+} \over \ds
{1\over\rho_+}- {1\over \rho_-}}}\,.\eeq
For a 1-shock one has  $\rho_+>\rho_-$, while for a 2-shock one has $\rho_+<\rho_-$.

The following  observation will be useful. Setting
\bel{sth} s~=~u_--u_+\,,\qquad\qquad \theta~=~{\rho_+\over\rho_-}~=~\left({h_+
\over h_-}\right)^{2\over \gamma-1},\eeq
from (\ref{RH2}) it follows that, for any shock wave,
\bel{RH4}
s~=~\frac{h_-}{B}\, \sqrt{A(1-\theta)(1-\theta^\gamma)\over\theta}\,,
\eeq
where $B$ is the constant at (\ref{xi_eq}).
%

\subsection{Small wave interactions}

Next, we review some well known results on wave interactions, for future use. 
Note that, when  a wave-front crosses a shock or a compression of the opposite family,
 the density $\rho$ (and hence $h$ as well) increases.   On the other hand, the density
 along a wave-front decreases when it crosses a rarefaction of the opposite family.
 
For any pairwise interaction between two small (shock or rarefaction) waves, one has the following estimates (see \cite{Bbook, GL, Sm}).

\v
\begin{proposition}\label{prop1}
Call $\sigma', \sigma''$ the strengths of two interacting wave-fronts, and let 
$\sigma_1, \sigma_2$ be the strengths of the outgoing waves  of the first and second family, in the solution of the Riemann problem. 
Then there exists a constant $C_0$ (uniformly valid as the state of the system ranges over a bounded set in the $\rho$-$u$ plane, with $\rho$ bounded away from zero) such that
\begin{itemize}
\item If $\sigma'$ is a 1-wave and $\sigma''$ is a 2-wave, then
\bel{299}
|\sigma_1-\sigma'|+|\sigma_2-\sigma''|~\leq~C_0\,|\sigma'\sigma''|(|\sigma'|+|\sigma''|).\eeq
\item If both $\sigma'$ and $\sigma''$ belong to the first family, then
\bel{21}
|\sigma_1-(\sigma'+\sigma'')|+|\sigma_2|~\leq~C_0\,|\sigma'\sigma''|(|\sigma'|+|\sigma''|).\eeq
\item If both $\sigma'$ and $\sigma''$ belong to the second family, then
\bel{22}
|\sigma_1|+|\sigma_2-(\sigma'+\sigma'')|~\leq~C_0\,|\sigma'\sigma''|(|\sigma'|+|\sigma''|).\eeq
\end{itemize}
\end{proposition}

\subsection{A rarefaction or compression wave crosses a large shock.}
\label{sec3.2}

\begin{figure}[htbp]
\centering
  \includegraphics[scale=0.45]{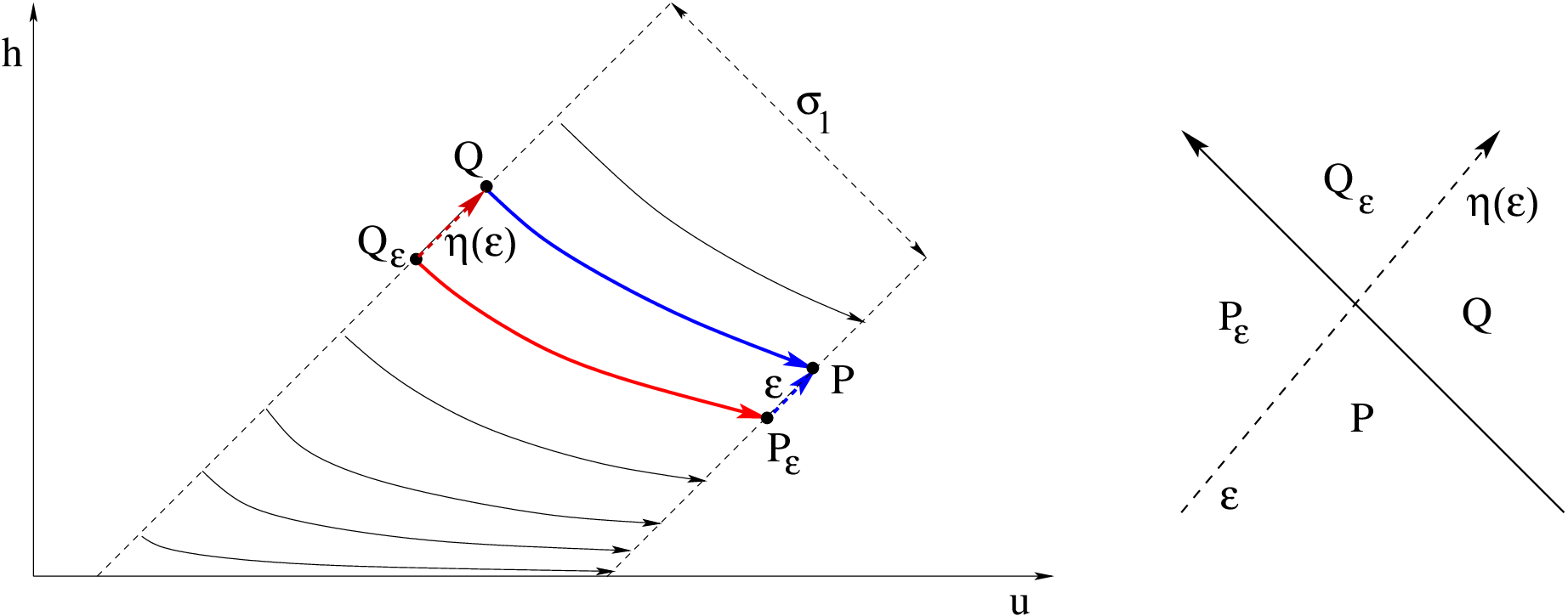}
    \caption{{\small A small 2-rarefaction crosses a 1-shock.
    In this configuration, the strength of the shock does not change, 
    while the strength of the rarefaction increases from $\ve$ to some value 
    $\eta(\ve)>\ve$.}}
\label{f:hyp207}
\end{figure}
To fix the ideas, consider a large 1-shock which
crosses a small 2-wave (compression or rarefaction)
of size $\sigma_2=\ve$.
We seek an estimate on the size of the outgoing waves, up to leading order.
As shown in Fig.~\ref{f:hyp207}, let
$$P~=~(u_-, \,h_-),\qquad\qquad Q~=~(u_+,\,h_+)
$$
be the left and right states across the large 1-shock before the interaction, and
let 
$$P_\ve~=~(u_--\ve, \,h_--\ve),\qquad\qquad Q_\ve~=~\Big(u_+-\eta(\ve),
\,h_+-\eta(\ve)\Big)
$$
be the left and right states across the 1-shock after the interaction.
Set
\bel{ff}s(\ve)~\doteq~(u_--\ve) - \bigl(u_+ - \eta(\ve)\bigr),\qquad
\theta(\ve)~\doteq
~\left({h_+ -\eta(\ve)
\over h_--\ve}\right)^{2\over \gamma-1}.\eeq
By (\ref{RH4}), replacing $P$ with $P_\ve$   we can write 
\bel{RH7}
s(\ve)~=~\frac{h_--\ve}{B}\, \sqrt{\psi(\theta(\ve))},\qquad 
\hbox{with}\qquad \psi(\theta)~\doteq~
{A(\theta-1)(\theta^\gamma-1)\over \theta}\,.\eeq
Hence
\beq
\eta(\ve)-\ve~=~s(\ve) - (u_- -u_+)~=~\frac{h_--\ve}{B}\sqrt{\psi(\theta(\ve))}- (u_- -u_+)\,.
\eeq
Differentiating w.r.t.~$\ve$ we obtain
\bel{der}\bega{rl}
\eta'(\ve)-1&\ds=~-\frac1B\sqrt{\psi(\theta(\ve))} + {h_--\ve \over 2B\sqrt{\psi(\theta(\ve))}}\cdot \psi'(\theta(\ve))\theta'(\ve),
\enda\eeq
\bel{26}
\theta'(\ve)~=~\frac2{\gamma-1}[\theta(\ve)]^{(3-\gamma)/2}\cdot
{-\eta'(\ve)(h_--\ve)+(h_+ -\eta(\ve)) \over (h_--\ve)^2}
\,.
\eeq
By (\ref{ff}),  at $\ve=0$ the above expression reduces to
\bel{266} \theta'~=~\frac2{\gamma-1}\theta^{(3-\gamma)/2}\cdot
{-\eta' +\theta^{(\gamma-1)/2} \over h_-}
\,.
\eeq
Using  (\ref{26}) to compute  the right hand side of (\ref{der}), when $\ve=0$ 
we find 
\bel{28}
\eta'-1 ~=~
-\frac1B\sqrt{\psi(\theta)} + {1\over 2B\sqrt{\psi(\theta)}}\cdot \psi'(\theta)
\cdot \frac2{\gamma-1}\theta^{(3-\gamma)/2}\cdot
\Big(-\eta' +\theta^{(\gamma-1)/2} \Big)
\,,\eeq
Solving for $\eta'$, we finally obtain
\bel{eta'}
\eta'~=~{2B\sqrt\psi-2\psi+\frac2{\gamma-1}\theta \psi'\over 2B\sqrt \psi+\frac2{\gamma-1}\theta^{(3-\gamma)/2}\psi'}~\doteq ~a(\theta).
\eeq

We observe that $\eta'$ is the factor by which an infinitesimal 2-wave
(either a compression or a rarefaction) is amplified when it crosses the 1-shock.
According to (\ref{eta'}), this ratio depends only on $\theta$. In particular, 
as $\theta$ remains bounded,  the above amplification coefficient is a bounded number.
To compute the amplification of an arbitrary  rarefaction or compression wave 
which crosses a large shock of  the opposite family, we can simply integrate
(\ref{eta'}) and  obtain
\bel{eee}\eta(\bar \ve)~=~\int_0^{\bar \ve} a(\theta(\ve))\, d\ve\,.\eeq


By a direct calculation we now prove $a(\theta)>1$.   In other words, as a
compression or rarefaction wave crosses a shock of the opposite family,
its strength always increases.  
In view of (\ref{eta'}), this will be a consequence of  the two inequalities 
\bel{41}
\left(2B\sqrt{\psi}-2\psi+\frac2{\gamma-1}\theta\psi'\right)-\left(2B\sqrt{\psi} +\frac2{\gamma-1}\theta^{3-\gamma\over 2}\psi'\right)~>~0,
\eeq
\bel{42}
2B\sqrt{\psi} +\frac2{\gamma-1}\theta^{3-\gamma\over 2}\psi'~>~0.
\eeq
We begin by observing that
 $\theta= h_+/h_->1$, and hence
\bel{1stder}
\psi(\theta)~\doteq~
{A(\theta-1)(\theta^\gamma-1)\over \theta}~>~0\,,
\qquad\qquad \psi'(\theta)~=~A\theta^{-2}[-1+(1-\gamma)\theta^\gamma+\gamma\theta^{\gamma+1}]~>~0,\eeq 
proving (\ref{42}).
 Moreover, one has
\bel{43}
\bega{l}
\ds\left(2B\sqrt{\psi}-2\psi+\frac2{\gamma-1}\theta\psi'\right)-
\left(2B\sqrt{\psi} +\frac2{\gamma-1}\theta^{3-\gamma\over 2}\psi'\right)\cr\cr
\qquad =~\ds-2\psi+\frac2{\gamma-1}\theta\psi'-\frac2{\gamma-1}\theta^{3-\gamma\over 2}\psi'\cr\cr
\qquad =~\ds-2A\bigl(\theta^{-1}+\theta^\gamma-\theta^{\gamma-1}-1\bigr)+\frac2{\gamma-1}\Big(-A\theta^{-1}+A(1-\gamma)\theta^{\gamma-1}+A\gamma\theta^\gamma\Big)(1-\theta^{1-\gamma\over2})\cr\cr
\qquad =~\ds \frac{2A}{\gamma-1}\bigl(1-\theta^{-{\gamma+1\over
2}}\bigr)
\Big[-\gamma\theta^{\gamma-1\over2}(\theta-1)
+(\theta^\gamma-1)\Big]~>~0.
\enda
\eeq
Indeed, to see that the last factor on the right hand side of 
(\ref{43}) is positive for $\theta>1$, we set
$$f(\theta)~=~-\gamma\theta^{\frac{\gamma-1}2}(\theta-1)+\theta^\gamma-1~=~-\gamma\theta^{\frac{\gamma+1}2}+\gamma\theta^{\frac{\gamma-1}2}+\theta^{\gamma}-1\,.$$
Then 
\[f(1)=0,\qquad\qquad 
f'(\theta)~=~\gamma\theta^{\frac{\gamma-3}2}\left[
-\frac{\gamma+1}2\theta+\frac{\gamma-1}2+\theta^{\frac{\gamma+1}2}\right]~=~
\gamma\theta^{\frac{\gamma-3}2}\cdot g(\theta),
\]
where
$$g(\theta)~=~-\frac{\gamma+1}2\theta+\frac{\gamma-1}2+\theta^{\frac{\gamma+1}2}\,,$$  
$$g(1)=0,\qquad\qquad g'(\theta)~=~-\frac{\gamma+1}2+\frac{\gamma+1}2\theta^{\frac{\gamma-1}2}~>~0.$$
For $\theta>1$ we thus have
$g(\theta)>0$, hence  $f'(\theta)>0$ and  $f(\theta)>0$.
This completes the proof that $a(\theta)>1$ for $\theta>1$.

With reference to Figure~\ref{f:hyp207}, the  inequality $|\eta(\ve)|>|\ve|$ 
implies that the $h$-components of the states $P,Q,P_\ve,Q_\ve$ satisfy
\bel{result1} |h_Q-h_{Q'}|~>~|h_P-h_{P'}|\,.\eeq

\subsection{A small shock crosses a large shock}

\begin{figure}[htbp]
\centering
  \includegraphics[scale=0.45]{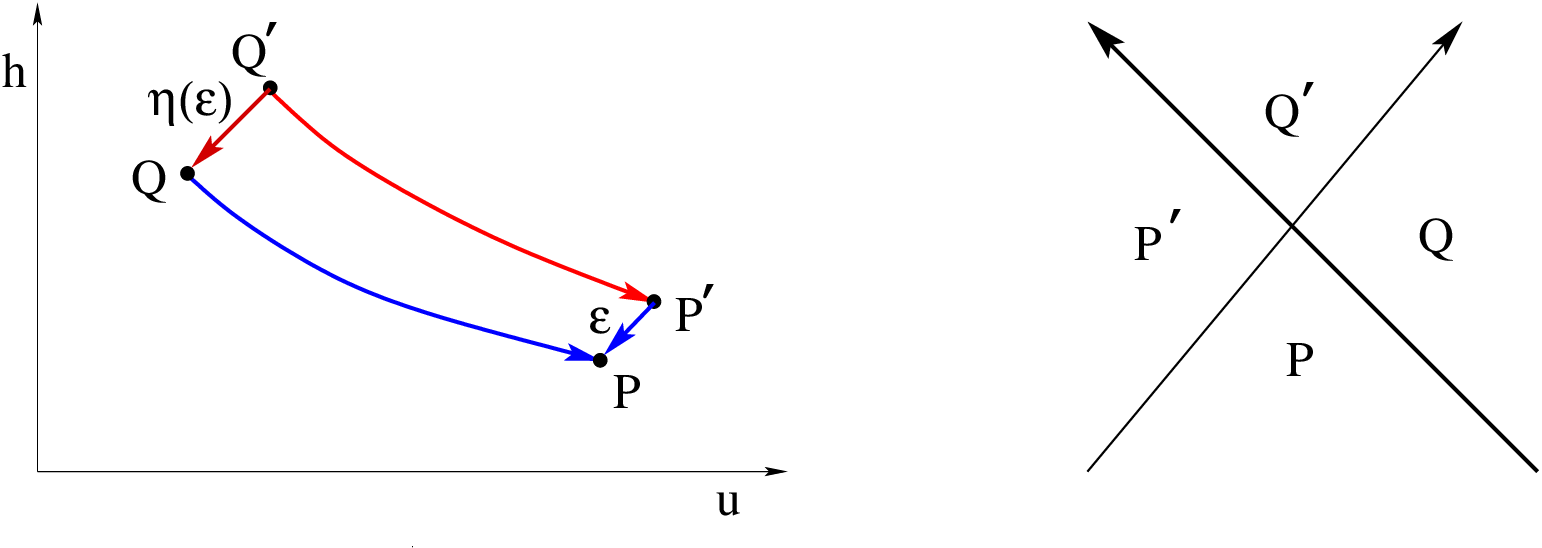}
    \caption{{\small A small 2-shock crosses a large 1-shock.}}
\label{f:hyp206}
\end{figure}

Let $\ve<0$ be the  signed strength of the small shock.   
Since shock and rarefaction curves
have a second-order tangency
\cite{Bbook, Sm}, with the notation used in Fig.~\ref{f:hyp206}
we have
$$P~=~(u_-, \,h_-),\qquad\qquad Q~=~(u_+,\,h_+),
$$
$$P'~=~(u_--\ve+o(\ve^2)~, ~\,h_--\ve+o(\ve^2)),\qquad\qquad Q'~=~\Big(u_+-\eta(\ve)+o(\ve^2),
\,h_+-\eta(\ve)+o(\ve^2)\Big).
$$
Here and in the sequel, the Landau notation $o(\ve^2)$
denotes an infinitesimal of higher order w.r.t.~$\ve^2$.
Computing the derivative $\eta'$ at $\ve=0$, we
thus recover exactly the same expression as in (\ref{28}).

Because of the second order tangency condition, 
the change in the strength of the big shock
will also be of order $o(\ve^2)$. 
\v
We conclude with an estimate which will be used later.
\begin{lemma}\label{lemma0}
Fix $0<a<b$ and consider the Riemann problem 
determined by the interaction of a large 1-shock 
with another wave front of strength $|\sigma|$. Assume that
the left, middle, and right states remain in the region  
where $a\leq h(v)\leq b$.  Call $|\sigma'|$
the strength of the outgoing 2-wave generated by the interaction. Then there exists a constant $C_\gamma$ depending only on $a,b$ such that:
\begi
\item[(i)] If the second wave impinges on the 1-shock from the left, then strength of the outgoing 2-wave satisfies
$$|\sigma'|~\leq~C_\gamma\,|\sigma|.$$
\item[(ii)] If the second wave is a 1-shock, or a small  1-compression  
or 1-rarefaction, 
impinging on the 1-shock from the right, then strength of the outgoing 2-wave satisfies
$$|\sigma'|~\leq~|\sigma|.$$
\endi
\end{lemma}

{\bf Proof.} Part (i) is an immediate consequence of Proposition~\ref{prop1}. 

Toward a proof of (ii), we first recall a basic property of shock curves 
for the p-system. 
Fix a left state $(u_-, h_-)$, and consider the curve of all 
points $(u_+, h_+)$ which can be connected to $(u_-, h_-)$ by a 1-shock.
Writing $h_= = h_+(u_+)$, the slope of this curve satisfies
\bel{sslope} -1~\leq~{dh_+\over d u_+} <~0.\eeq
Indeed, this inequality  is established within the proof of Lemma 3 in  \cite{NS}.
It is also found in Section~3 of \cite{CJ}.

Next, assume  that the left and right states across the large 1-shock are
$$P~=~(u_P, \,h_P),\qquad\qquad Q~=~(u_Q,\,h_Q),
$$
Two cases will be considered.

\begin{figure}[htbp]
\centering
  \includegraphics[scale=0.65]{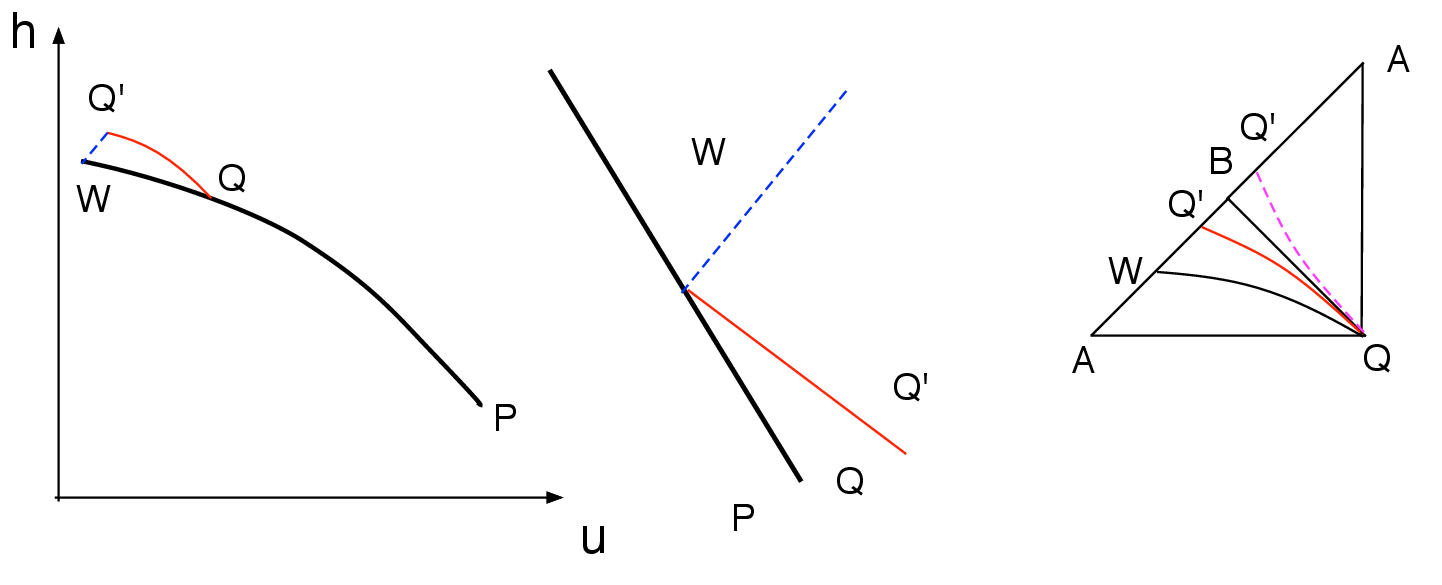}
    \caption{{\small A small 1-shock impinges a large 1-shock.  }}
\label{f:p207}
\end{figure}

CASE 1:  The impinging 1-wave is a shock, 
connecting the states $Q$, $Q'$.
As shown in Fig.~\ref{f:p207}, the outgoing 2-rarefaction connects the 
states $W$, $Q'$, where $W=(u_W, h_W)$ is the unique state along 
the 1-shock curve through $P$ such that 
$$u_{Q'}- h_{Q'} ~=~u_W- h_W.$$
In this case, call $A$ and $B$ the points in the $u$-$h$ plane such that 
$$u_{Q'}- h_{Q'} ~=~u_A- h_A~=~u_B-h_B,$$
$$h_A ~=~h_Q\,,\qquad\qquad u_B+h_B~=~u_Q + h_Q\,.$$
We then have
 $$\bega{l}\hbox{[strength of the outgoing 2-rarefaction]}~=~\Big|(u_{Q'}+ h_{Q'}) - 
(u_W+h_W)\Big|\\[4mm]
\qquad \leq~
(u_B+ h_B) - 
(u_A+h_A)
=~(u_Q-h_Q) - (u_B-h_B)\\[4mm]\qquad =~(u_Q-h_Q) - (u_{Q'}-h_{Q'})
~=~\hbox{[strength of the incoming 1-shock].}
\enda
$$

CASE 2:
The incoming 1-wave is a compression. In this case we have 
$Q'=B$, and the previous inequalities remain valid.

\begin{figure}[htbp]
\centering
  \includegraphics[scale=0.65]{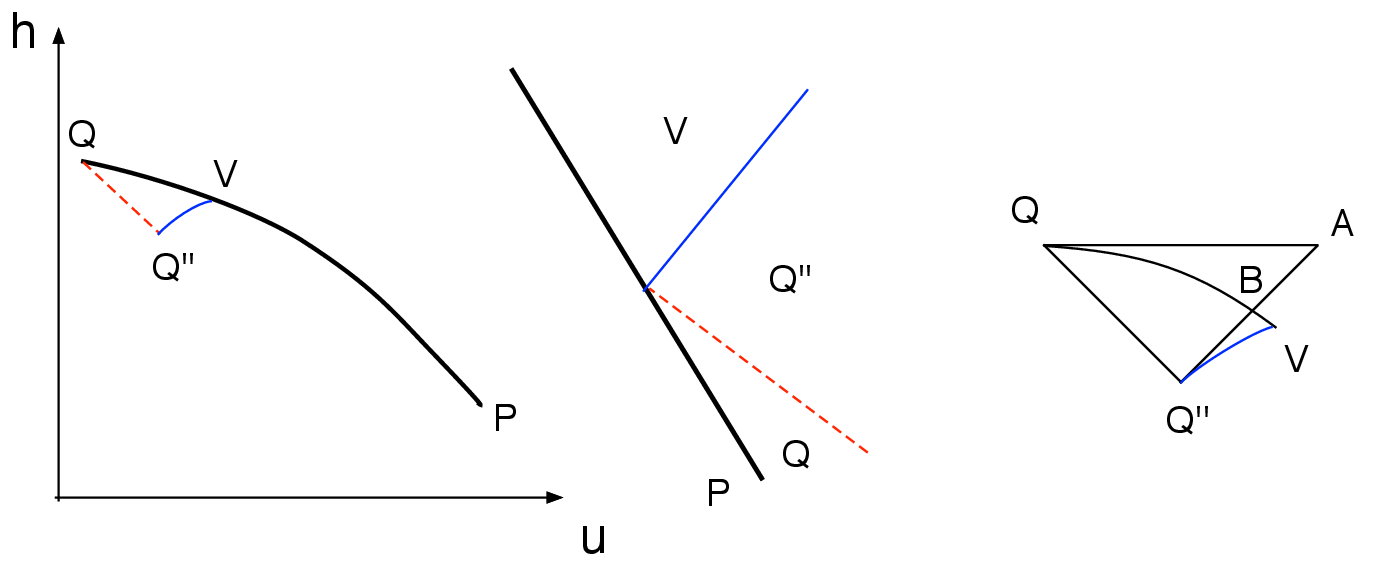}
    \caption{{\small A small 1-rarefaction impinges a large 1-shock. }}
\label{f:p208}
\end{figure}

CASE 3: The impinging 1-wave is a small rarefaction, connecting $Q$ with a right state $Q''$, of strength $$\ve~=~(u_{Q''}- h_{Q''}) - (u_Q-h_Q).$$
In this case the interaction 
produces an outgoing 2-shock, connecting the states $V$ and $Q''$.
Here $V$ is the state at the intersection of the 1-shock curve through $P$
and the 2-shock curve through $Q''$.
As shown in Fig.~\ref{f:p208}, call $A= (u_A, h_A)$ the point such hat
$$u_A-h_A ~=~u_{Q''} - h_{Q''}\,,\qquad\qquad h_A~=~h_Q\,,$$
and let $B$ be the point at the intersection of the 1-shock curve
through $P$ and the segment $AQ''$.
Recalling that the 2-shock curve through $Q''$ has a second order tangency 
with the segment $AQ''$, and using the inequality
$$-1~\leq~{dh+\over du_+}~\leq~-\delta_0 ~<~0$$
for some $\delta_0>0$ sufficiently small,
we obtain
 $$\bega{l}\hbox{[strength of the outgoing 2-shock]}~=~(u_V+ h_V) - 
(u_{Q''}+h_{Q''})\\[4mm]
\qquad =~(u_B+ h_B) - 
(u_{Q''}+h_{Q''})+ \O(\ve^3)\\[4mm]\ds
\qquad \leq~\left(1-{\delta_0\over 2}\right)
\Big[(u_{Q''}- h_{Q''}) - (u_Q-h_Q)\Big] +\O(\ve^3)\\[4mm]\qquad 
\ds=~\left(1-{\delta_0\over 2}\right)\ve +\O(\ve^3)~<~\ve~=~\hbox{[strength of the incoming 1-rarefaction].}
\enda
$$
Together, the above three cases prove part  (ii) of Lemma \ref{lemma0}.\endproof


\subsection{Wave measures.} 
\label{s:wm} Let now $x\mapsto (v(x),u(x))$ be any 
profile with bounded variation. As in Chapter 10 
of \cite{Bbook}, we can define 
the signed measures $\mu_1, \mu_2$ describing strength of waves.
Namely, for $i=1,2$, 
\begi
\item The atomic part of $\mu_i$ is supported on the countable set of
points where $v$ or $u$ have a jump.   If $\bar x$ is one such point, then
$\mu_i(\{\bar x\})$ is the signed strength of the $i$-th wave in the solution of the Riemann problem with left and right data
$$(v,u)(\bar x -),\qquad (v,u)(\bar x +).$$
\item The continuous part of $\mu_i$ is defined as the continuous part 
of the distributional derivative of the scalar function $x\mapsto w_i(x)$.
Since $w_i$ has bounded variation, this is a bounded measure.
\endi
These measures can be decomposed into a positive and a negative part, so that
$$\mu_i~=~\mu_i^+ - \mu_i^-\,,\qquad\qquad 
|\mu_i|~=~\mu_i^+ + \mu_i^-\,.$$
Notice that:
\begi\item $\mu_i^+$ accounts for $i$-rarefaction waves.
\item The continuous part of 
$\mu_i^-$ accounts for $i$-compression waves
\item The atomic part of 
$\mu_i^-$ accounts for $i$-shocks.
\endi

As shown in \cite{Bbook}, Glimm's functionals (originally 
defined for piecewise constant functions) can be 
extended to
arbitrary BV functions.   
The {\bf total strength of waves} is defined as
\bel{V}
V~\doteq~\sum_{i=1,2}|\mu_i|(\R),\eeq
while the {\bf interaction potential} is
\bel{Q}Q~\doteq~\int_{x<y}d|\mu_2|(x)\, d|\mu_1|(y)+
\sum_{i=1,2} |\mu_i|(\R)\cdot\mu_i^-(\R)\,.\eeq
Notice that (\ref{Q}) accounts for the product of
strengths of all couples of approaching waves.
We recall that two waves of the same family are
{\it approaching} if at least one of them is a compression or a shock.

Next, consider a solution $(v,u)$ of (\ref{1}) defined for 
$t\in [t_0,t_1]$ and
let $V(t)$, $Q(t)$ be the total strength of waves and the 
wave interaction potential at time $t$. 
As shown in Chapter 10 of \cite{Bbook}, these functionals
satisfy the same estimates valid for Glimm or front-tracking
approximations.   In particular,
from the interaction estimates in Proposition~\ref{prop1}
it follows 
\begin{lemma}\label{lemma9}
For any given $K_0$ and $a,\ve_0>0$, 
there exists $\delta_0>0$ such that
the following holds.  Assume that
\begi
\item[(i)] the density remains bounded away from zero: 
$h(v(t,x))\geq a$ for all $t\in [t_0, t_1]$, $x\in\R$,
\item[(ii)] the total strength of waves in the initial data satisfies
$V(t_0)\leq K_0$, and 
\item[(iii)]
the solution $(v,u)$ does not 
contain any shock of strength $>\delta_0$.
\endi
Then the function
\bel{GF}
t~\mapsto~V(t)+\ve_0 Q(t),\qquad\qquad t\in [t_0, t_1]\eeq
is non-increasing.
\end{lemma}
\v
 Given a BV solution 
$U=(v,u)$ of (\ref{1}), 
we denote by
\bel{l12}\lambda_1(t,x)~\doteq~-c(v(t,x)),
\qquad\qquad\lambda_2(t,x)~\doteq~c(v(t,x)),\eeq
the two wave speeds at the point $(t,x)$, as in
(\ref{speed}).
Following \cite{Daf}, for $i=1,2$, by a {\it generalized $i$-characteristic} 
we mean an absolutely continuous curve
$x=x(t)$ such that
\bel{48}\dot x (t)~\in~ \bigl[ \lambda_i(t,x+),~
\lambda_i(t,x-)\bigr]\eeq
for a.e.~$t$.

For a given terminal point $\bar x$ we shall consider the
{\it minimal $i$-characteristic} through $\bar x$, defined as
$$\xi(t,\bar x)~\doteq~\min\big\{ x(t)\,; ~~x~\hbox{ is an $i$-characteristic},
~~x(t_1)=\bar x\big\}.$$
As proved in \cite{Daf}, the curve $\xi(\cdot,\bar x)$ is itself an $i$-characteristic. Indeed, for a.e.~$t$ the functions
$w_1,w_2$ and hence also the wave speed
$\lambda_i$ are continuous
at $(t, \xi(t))$. Therefore we can simply write
$$\dot \xi(t)~=~\lambda_i(t, \xi(t)).$$
\v
In addition to the wave measures $\mu_i$, one can also 
introduce a scalar, positive measure $\mu^{int}$ on the 
$t$-$x$ plane bounding the amount of interaction, and hence
the production of new waves.  More precisely, 
let $(U_\nu)_{\nu\geq 1}$ be a sequence of piecewise constant
front-tracking solutions, converging to the exact BV 
solution  $U=(v,u)$.    

For each $\nu\geq 1$ we can also construct 
a purely atomic measure 
$\mu^{int}_\nu$ by setting
\bel{mint}\mu^{int}_\nu\bigl(\{\ov P\}\bigr)~=~|\sigma'\,\sigma''|\eeq
for every point $\ov P=(\bar t, \bar x)$ at which 
two incoming fronts interact, with strengths $\sigma', \sigma''$
respectively.
By taking a suitable subsequence, we can achieve the weak convergence of measures
\bel{mint2}\mu_\nu^{int}~\to~\mu^{int}\eeq
for some positive measure $\mu^{int}$, which we call
a {\bf measure of wave interaction} for the solution $U$.
Taking the limit of front tracking approximations
one obtains a useful important
property of this measure, namely:

\begin{lemma}\label{lemma5}  For $i\in \{1,2\}$, let 
$t\mapsto \xi(t)$ and $t\mapsto \tilde\xi (t)$ be two minimal 
i-characteristics, with $ \xi (t)\leq  \tilde\xi (t)$ for $t\in [t_0, t_1]$.
Then one has the estimate
\bel{midif}
\mu_i^\pm\Big(\bigl[ \xi(t_1) , \,\tilde\xi(t_1)\bigr[\Big)~\leq~
\mu_i^\pm\Big(\bigl[\xi(t_0) , \,\tilde\xi(t_0)\bigr[\Big)+
C\cdot\mu^{int}(\Omega),\eeq
where
\bel{Omdef}\Omega~\doteq~
\Big\{ (t,x)\,;~~t\in [t_0, t_1]\,,\,
x\in [\xi(t), \,\tilde\xi(t)[\,\Big\}.\eeq
\end{lemma}
In other words, the total amount of (positive or negative) 
$i$-waves at time $t_1$
contained in the interval $[\xi(t_1), \,\tilde\xi(t_1)[$ can be estimated in terms of 
the ``old $i$-waves" (positive or negative, respectively) present at time $t_0$
inside the interval $[\xi(t_0) , \,\tilde\xi(t_0)[\,$, plus
some ``new waves" generated by wave interactions occurring
inside the domain $\Omega$ enclosed between the two characteristics.  The total strength of these new waves
can be bounded in terms of the interaction measure $\mu^{int}$.
\v
A precise value for the constant $C$ in
(\ref{midif}) can be determined using the interaction estimates 
(\ref{299})--(\ref{22}). In particular, by taking the limit of front tracking approximations, one obtains

\begin{lemma}\label{lemma7} 
For any given constants $K_0,a$ and
$\ve_0>0$, one can find $\delta_0>0$ such that, under the 
assumptions (i)--(iii) of Lemma~\ref{lemma9}, one has 
\begi
\item[(i)] the total amount of interaction satisfies
\bel{tint}\mu^{int}\bigl([t_0, t_1]\times\R\bigr)~\leq~2K_0^2\,,
\eeq
\item[(ii)] the estimate (\ref{midif})
 holds with $C=\ve_0$.
\endi
\end{lemma}

\v
\section{Decay of positive waves.}
\label{sec:22}
\setcounter{equation}{0}
Differentiating (\ref{w12}) and writing the wave speed as 
$c=c(w_2-w_1)$
with $c$ as in (\ref{cdef}), one obtains
\bel{wx}
\left\{\bega{rl} w_{1,xt} - c w_{1,xx}&
=~ -c' w^2_{1,x}+ c' w_{1,x} w_{2,x}
\,,\\[3mm]
w_{2,xt} + c w_{2,xx}&=~c'w_{1,x}  w_{2,x}- c' w^2_{2,x}
\,.\enda \right. 
\eeq
The above system would be easy to integrate if we did not have the 
mixed terms $w_{1,x}w_{2,x}$.  To get rid of these terms, we first multiply 
each equation in (\ref{w12}) by
a  function $\phi=\phi(w_2-w_1)$ and then differentiate.  
For example, the second equation yields
\bel{90}\bega{l}[\phi\, w_{2,x}]_t + c [\phi \,w_{2,x}]_x\\[3mm]
\qquad =~\phi\,\Big[c' w_{1,x}w_{2,x}- c' w^2_{2,x}\Big]
+ \phi'\,(w_{2,t}-w_{1,t})w_{2,x}
+ c\,\phi'(w_{2,x}-w_{1,x})w_{2,x}
\\[3mm]
\qquad =~-\phi\,\Big[ -c' w_{1,x}w_{2,x}+c_{w_2} w^2_{2,x}\Big]
- 2c \,\phi'\,w_{1,x}w_{2,x}\\[3mm]
 \qquad=~- \phi\, c' w^2_{2,x} \,,\enda
\eeq
provided that 
\bel{91}
{\phi'\over\phi} ~=~{c'\over 2 c}~=~\frac12\,\frac{\gamma+1}{\gamma-1}\,\frac1{w_2-w_1}\,.\eeq
Computing an explicit solution of (\ref{91}) we find
\bel{eqphi}
\phi(w_2-w_1)~=~(w_2-w_1)^{\frac{\gamma+1}{2\gamma-2}}~=~
(2h)^{\frac{\gamma+1}{2\gamma-2}}\,.
\eeq
 In the end, this yields
a decay estimate along any  2-characteristic $t\mapsto x(t)$.
\bel{dec1}
{d\over dt} (\phi w_{2,x})\bigl(t, x(t)\bigr)~=~-\phi \,c_{w_2}\, w_{2,x}^2~\leq~-C_2 w_{2,x}^2\,,\eeq
for some constant $C_2>0$ depending only on the 
upper and lower bounds for the density.  Of course, an entirely similar estimate
holds for 1-rarefactions.
\v
Next, let $\xi_1(t)<\xi_2(t)$ be two 2-characteristics. 
Calling $c=c(w_2-w_1)$ the characteristics speed as a function 
of the Riemann coordinates, we  have
\bel{ch2}{d\over dt}(\xi_2(t)-\xi_2(t))~=~\dot \xi_2(t)-
\dot \xi_1(t)~=~\int_{\xi_1(t)}^{\xi_2(t)} c'\,
(w_{2,x}-w_{1,x})\, dx\,.\eeq
Notice that the above identity  involves
also the 1-waves inside the interval $[\xi_1(t),\, \xi_2(t)]$.
We seek an equivalent way to express the distance between
two characteristics, which does not involve the 
contribution of intermediate 1-waves.  Toward this goal, 
consider the integral
\bel{wd}
Z(t)~\doteq~\int_{\xi_1(t)}^{\xi_2(t)}
\vp\, dx\,,
\eeq
where $\vp=\vp(w_2-w_1)=(w_2-w_1)^{-\frac{\gamma+1}{2\gamma-2}}$ which satisfies $c'\vp=-2c\vp'$. 
We
compute
\bel{ch3}\bega{rl}\ds{d\over dt}Z(t)&=~\ds
\dot\xi_2\vp(\xi_2)-\dot\xi_1\vp(\xi_1)+
\int_{\xi_1}^{\xi_2}\vp'\,(w_{2,t}-w_{1,t})\, dx\\[4mm]
&=~\ds c(\xi_2)\vp(\xi_2)- c(\xi_1)\vp(\xi_1)
-\int_{\xi_1}^{\xi_2}\vp'\, c\, (w_{2,x}+w_{1,x})\, dx
\\[4mm]
&=~\ds\int_{\xi_1}^{\xi_2}\Big[ c'(w_{2,x}-w_{1,x}) \vp
+ c\,\vp' (w_{2,x}-w_{1,x})\Big]\, dx
-\int_{\xi_1}^{\xi_2}\vp'\, c\, (w_{2,x}+w_{1,x})\, dx
\\[4mm]
&=~\ds\int_{\xi_1}^{\xi_2}\Big[-2c\vp' (w_{2,x}-w_{1,x})
+ c\,\vp' (w_{2,x}-w_{1,x})\Big]\, dx
-\int_{\xi_1}^{\xi_2}\vp'\, c\, (w_{2,x}+w_{1,x})\, dx\\[4mm]
&=~\ds\int_{\xi_1}^{\xi_2} c'\vp\, w_{2,x}\, dx\,.
\enda
\eeq
Notice that the last two equalities were obtained using the identity
$c'\vp=-2c\vp'$, which produces a cancellation of all terms involving 
1-waves.

\v
\section{A periodic interaction pattern}
\label{sec:2}
\setcounter{equation}{0}

As a preliminary to the blow-up example, in this section we construct a piecewise constant 
approximate solution with a periodic interaction pattern.

\begin{figure}[ht]
\centerline{\hbox{\includegraphics[width=15cm]{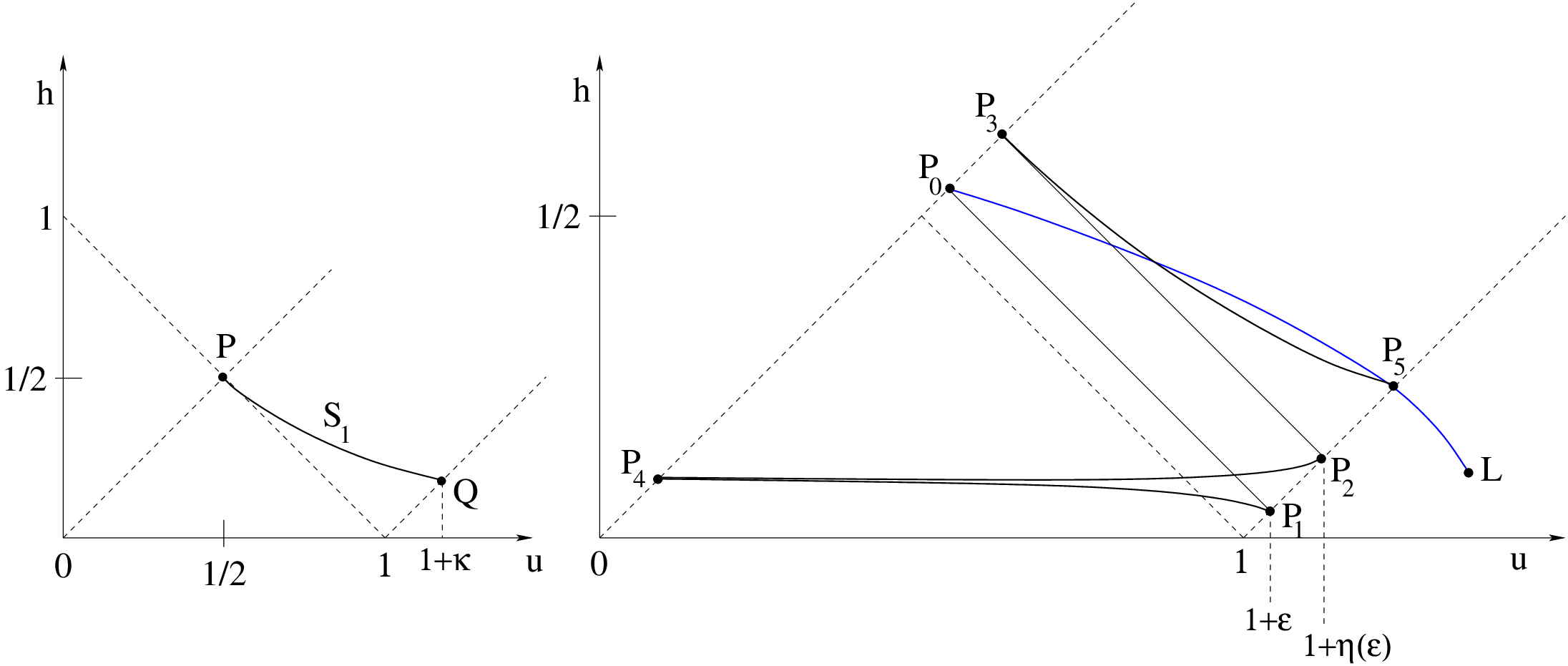}}}
\caption{\small Left: the 1-shock curve through $(\frac12,\frac12)$. Right: the 
various points $P_i$ considered in Lemma~1. }
\label{f:p68}
\end{figure}

In the following (see Fig.~\ref{f:p68}), we consider points $P_i=(u_i, h_i)$ along the two lines
$$\gamma_0~=~\{(u,h);\qquad h>0,~~u-h=0\},\qquad\qquad \gamma_1~=~
\{(u,h);\qquad
h>0,~~u-h=1\}.$$
\v
{\bf Lemma 2.} {\it 
There exists a point $P_0\in \gamma_0$ such that the following holds
(Fig.~\ref{f:p68}, right).  Consider the point 
$P_1=\Big(u_0+{1\over 2}, ~h_0-{1\over 2}\Big)
\in\gamma_1$.   Let $P_4\in \gamma_0$ be the point along the 1-shock
curve through $P_1$ and let $P_2\in \gamma_1$ be the point
along the 2-shock curve through $P_4$.
Finally, call $P_3=\Big(u_2-{1\over 2}, ~h_2+{1\over 2}\Big)
\in\gamma_0$ and let $P_5\in \gamma_1$ be the point along the 1-shock curve through $P_3$.
Then $h_5<h_0$.

As a consequence, there is a left state $L$ which can be connected to both 
$P_0$ and $P_5$ by 1-shocks.
}

{\bf Proof.}
First, consider the 1-shock through  
$P=\Big({1\over 2}, {1\over 2}\Big)$. 
This intersects the line   $\gamma_1$ 
at some point $Q$, say    with $h(Q)  =  \kappa>0$.

Next, for $\ve>0$ small consider the points 
$$P_0~=~\Big({1\over 2} +\ve,\, {1\over 2} +\ve\Big),
\qquad\quad    P_1~ =~ (1+ \ve,\, \ve).    $$
Starting from $P_1$, construct the corresponding point $P_4$ and
then  $P_2= (1+\eta(\ve), \,\eta(\ve))$.
By (\ref{eee}) and the boundedness of the 
amplification factor $a(\theta)$, as $\ve\to 0+$ we also have 
$\eta(\ve)\to 0$.

Finally, call $P_3 = (1+\eta(\ve),  1+\eta(\ve))$ and let 
$P_5\in \gamma_1$ be the point along the 1-shock curve through $P_3$.

By \eqref{result1}, 
it follows  
$$h_3 -   h_5~>~ h(P) - h(Q) ~=~   {1\over 2} - \kappa. $$ 
Therefore
$$\bega{l}h_0-h_5~=~(h_0-h_3)+(h_3-h_5)\\[4mm]\qquad 
=~(h_1-h_2)+(h_3-h_5)
>~\ve - \eta(\ve) +\Big({1\over 2}-\kappa\Big)   ~ >~ 0\enda$$   
for $\ve>0$  small enough.  We thus have
\bel{66}
h_5<h_0<h_3\,,\qquad\qquad  u_0<u_3< u_5\,.
\eeq
This proves the first statement in Lemma 2.

\begin{figure}[htbp]
\centering
  \includegraphics[scale=0.6]{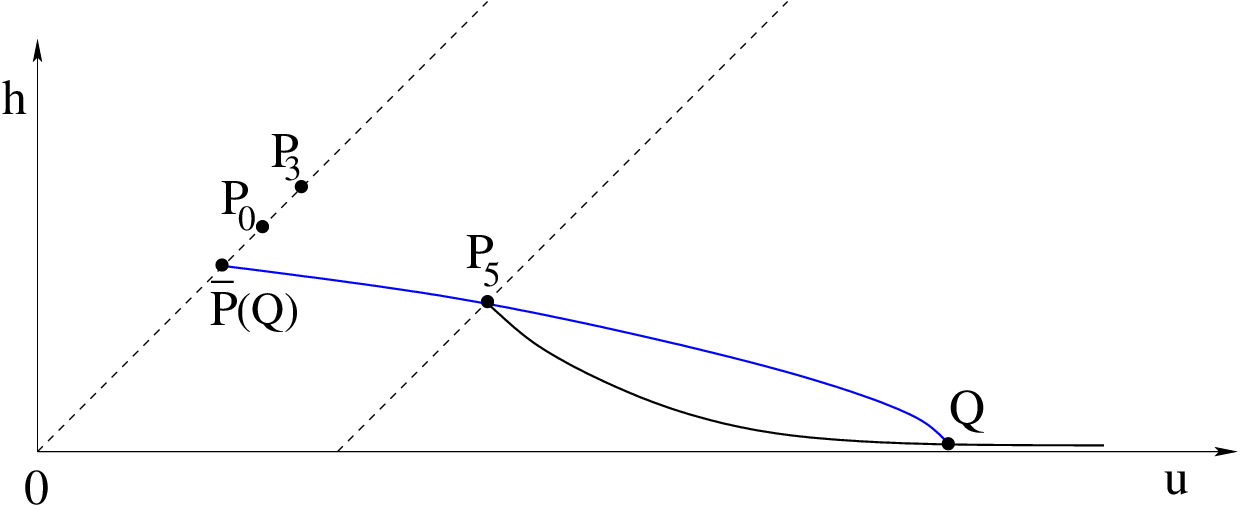}
    \caption{\small Constructing the left state $L$. }
\label{f:p60}
\end{figure}

It remains to prove that
there exists a left state $L= (u_L, \rho_L)$ which is connected to 
both $P_5$ and $P_0$ by a 1-shock. 
Referring to Fig.~\ref{f:p60} consider the 1-shock curve through $P_5$.
Let $Q= (u_Q, h_Q)$  be any point on this curve. Notice that, as $u_Q\to +\infty$, 
we have $h_Q\to 0$.

Next, let $\ov P(Q)= (\bar u, \bar h)$ be the point where  the 1-shock curve through $Q$   intersects the line $\gamma_0 = \{ h-u=0\}$.
Since this  shock curve is concave down, one has
$$\liminf_{u_Q\to +\infty} \,{ h_5-\bar h\over u_5-\bar u}~\geq~
\lim_{u_Q\to +\infty} \, {h_Q-h_5\over u_Q-u_5}~=~0.$$
Therefore, as $u_Q\to +\infty$ the point $\ov P(Q)$ approaches
the point $(u_5-1, h_5)$.  On the other hand, as $u_Q\to u_5$ one has 
$\ov P(Q)\to P_3$.   Since $h_0>h_5$, by continuity, there is some 
choice of $Q$ such that $\ov P(Q)= P_3$.   This completes the proof of Lemma 2.
\endproof

\begin{figure}
	\centering
	\includegraphics[scale=.6]{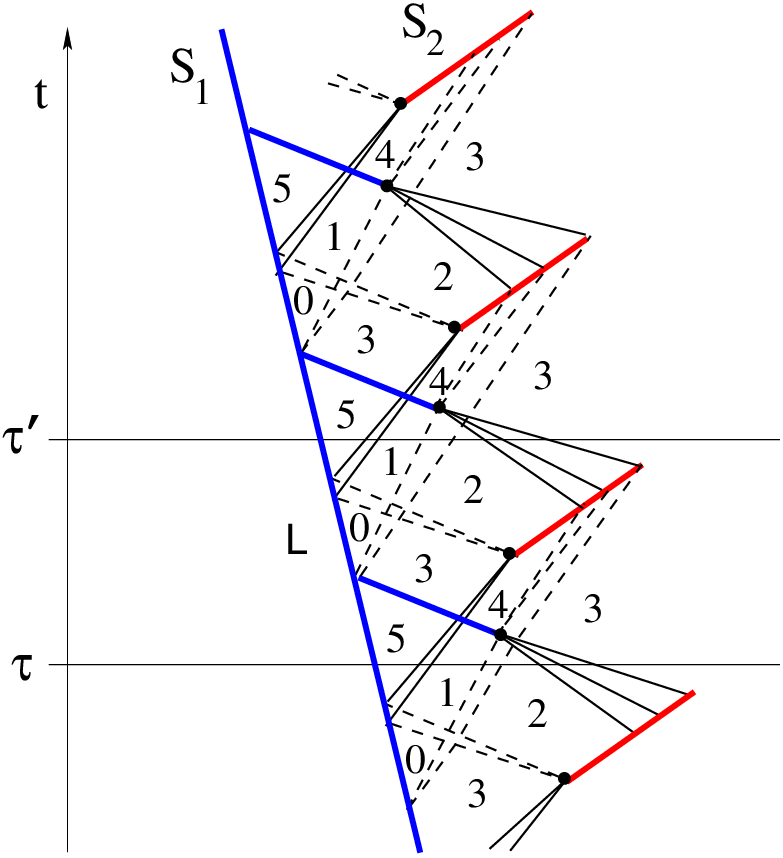}
	\caption{\small A periodic interaction pattern in the $(x,t)$-plane.
	Here  $0,1,\ldots, 5$ and $L$ refer to the 
	states $P_0,P_1,\ldots, P_5$ and to the left state 
	$L$ considered in Lemma~1.  The thick solid lines are shocks, 
	the thin solid lines represent compressions, 
	while the dashed lines are rarefaction fronts.
\label{f:p96} }
\end{figure}

Using Lemma~1 we now construct a front tracking solution
to the system (\ref{1}) with a periodic interaction pattern.
Referring to  Fig.~\ref{f:p96}, at time $t=\tau$ the piecewise constant 
solution $(u,h)(\tau,\cdot)$ takes the values
$L, P_5, P_2,P_3,P_4$.   As time increases, the following 
interactions
take place, one after the other.
\begi
\item[(i)] The 2-rarefaction $P_1P_2$ crosses the 1-compression $P_2P_3$. 
Afterwards, this 1-compression breaks into a 1-shock and a 2-rarefaction. 
In the end, the states $P_1$ and $P_3$ are connected 
by the 1-shock $P_1P_4$ followed by  the 2-rarefaction $P_4 P_3$.

\item[(ii)] The 2-compression $P_5P_1$ crosses the 1-shock
$P_1P_4$.   The Riemann problem is solved by the 1-shock $P_5P_3$
and the 2-compression $P_3P_4$.
\item[(iii)] The 1-shock $P_5P_3$ hits the 1-shock $L \,P_5$, generating 
the 1-shock $L \,P_0$ and the 2-rarefaction $P_0 P_3$.
\item[(iv)] The 2-compression $P_3P_4$ breaks into the 1-rarefaction 
$P_3P_2$ and the 2-shock $P_2 P_4$.
\item[(v)] The 2-rarefaction $P_0P_3$ crosses the 1-rarefaction
$P_3P_2$, producing the 1-rarefaction $P_0P_1$ and the 2-rarefaction 
$P_1P_2$.
\item[(vi)]
The 1-rarefaction $P_0 P_1$ hits the 1-shock $L\, P_0$,
producing the 1-shock $L\,P_5$ and the the 2-compression $P_5 P_1$.
\item[(vii)] The 2-shock $P_2P_4$ is canceled by 
the 2-rarefaction $P_4P_3$,
producing the  1-compression $P_2P_3$.
\endi
At time $t=\tau'$ we have reached the same configuration 
as at time $\tau$,
and the  periodic pattern can be continued.
\v
{\bf Remark.}  If the initial data had
small total variation, 
then the standard wave interaction estimates 
\cite{Bbook, G, Sm} 
would imply that 
the interaction potential approaches zero. As proved in
\cite{Liu}, the solution would
converge to the solution of the Riemann Problem with left and 
right data $(L, P_3)$.     In the present interaction pattern, however, this does not happen because wave strengths are large.
In particular, notice that  the 2-compression $P_5 P_1$
is greatly amplified when it crosses the large 1-shock 
$P_1 P_4$.

\section{An example with finite time blow-up of the total variation}
\label{sec:3}
\setcounter{equation}{0}

In this section we provide an affirmative answer to 
the question {\bf (Q)} considered in the Introduction.
Namely, we construct  a piecewise smooth approximate
solution of (\ref{1}) such that:
{\it 
\begi
\item[{\bf (C1)}] At each 
interaction, the strengths of outgoing waves is the same as in an exact solution.
\item[{\bf (C2)}] For some constant $C_0>0$, all
rarefaction waves satisfy a decay estimate of the form 
\bel{deci}{d\over dt} (\phi w_{i,x}) (t, x_i(t))~
\leq~-C_0 w_{i,x}^2\,,\qquad
\qquad i=1,2.\eeq
Here $\vp$ is the function at (\ref{eqphi}) and
$t\mapsto x_i(t)$ is any  $i$-characteristic.
\item[{\bf (C3)}] The density $\rho$ remains uniformly positive.
\item[{\bf (C4)}] The total variation blows up in finite time.
\endi
}

\begin{figure}
	\centering
	\includegraphics[scale=.45]{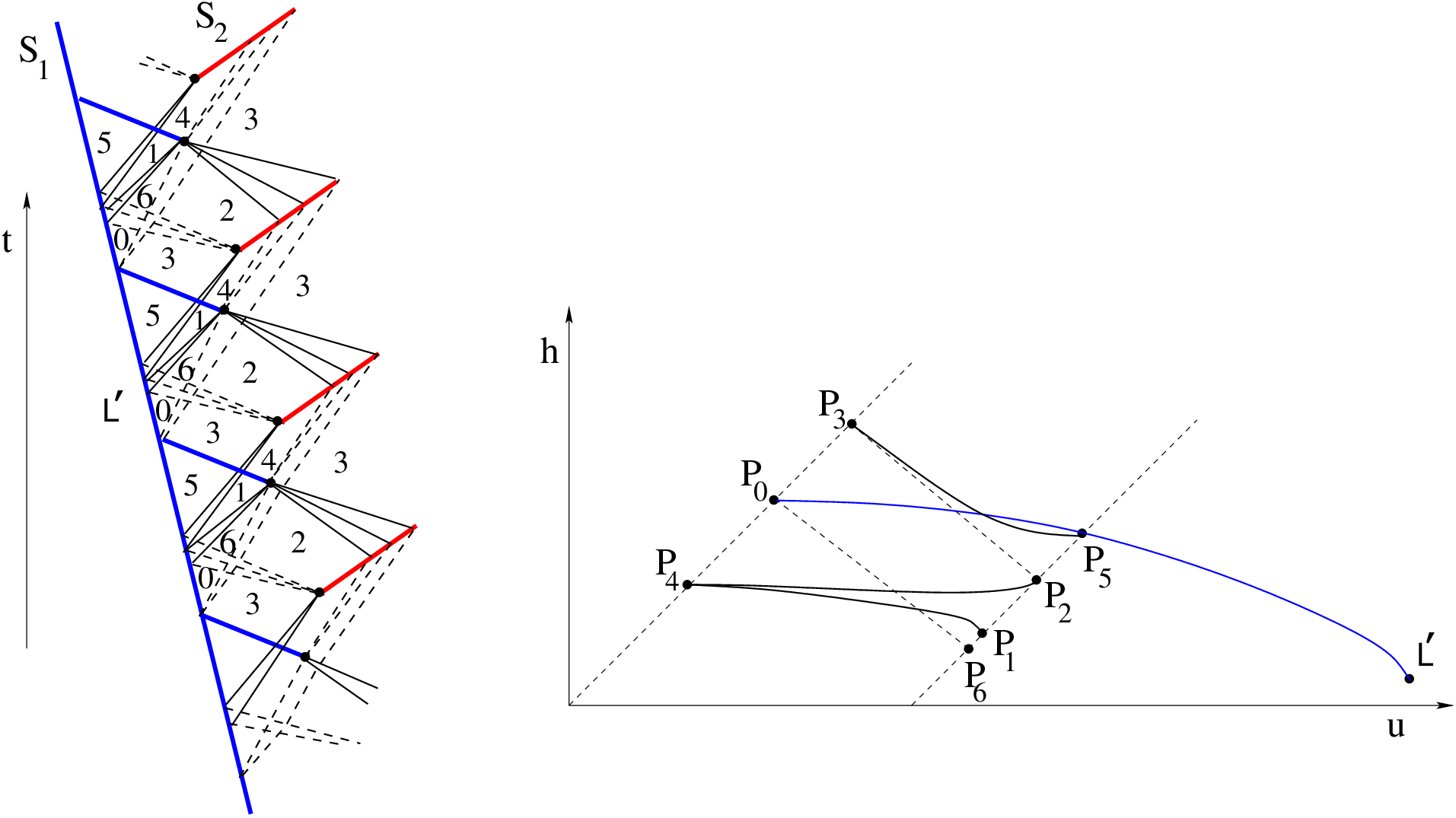}
	\caption{\small  Left: a modified periodic pattern.
	Compared with the pattern in  Fig.~\ref{f:p96},
	the compression wave joining $P_5$ to $P_1$ is now  
	split in two parts. This creates an additional state, which we call $P_6$.   Right: the location of the new states
	$P_0,\ldots, P_6$ and $L'$, in the $u$-$h$ plane.}
	\label{f:p90}.
\end{figure}

\subsection{Outline of the construction}

Then to construct  an approximate
solution of (\ref{1}) whose total variation blows  up in finite time,  the periodic 
pattern constructed in the previous section 
will be modified in two ways:
\v
\begi
\item[(i)] By slightly changing the wave speeds, 
the interaction 
pattern can be repeated on a sequence of shorter and shorter
time intervals $[\tau_{n-1}, \tau_n]$, with $\lim_{n\to\infty}\tau_n=T<\infty$. 
\item[(ii)] In the original pattern the 2-shock connecting $P_2$ with $P_4$
is entirely cancelled by the 2-rarefaction connecting $P_4$ with $P_3$.
We slightly change the speeds of these two waves so that they do not entirely 
cancel each other. More precisely, for every $n$ large enough, 
at the terminal time $T$ the solution will still contain a 2-shock and a 2-rarefaction, both of strength $\alpha\, n^{-1}$, connecting the  states $P_3^{(n+1)}$ and $P_3^{(n)}$. 
Here $\alpha$ is a fixed positive constant.
These are the remaining 
portions of the 1-shock and 1-rarefaction which are not completely cancelled 
by the $n$-th iteration of the basic pattern.
The total strength of all these waves is $\sum_{n\geq N} 2\alpha\,n^{-1}
= +\infty$, providing the blow-up of the total variation as $t\to T-$. 
\endi

Because of (ii), it is clear that the intermediate
states $P_0,\ldots P_5$ generated by this interaction pattern
can no longer repeat cyclically, but will 
slightly change after each round of interactions.  
Still, as $t\to T$, a periodic interaction 
pattern will be approached.   
The construction of 
the approximate solution will be achieved in the next three steps.
\v
\subsection{A perturbed periodic interaction pattern.}
\begin{figure}
	\centering
	\includegraphics[scale=.4]{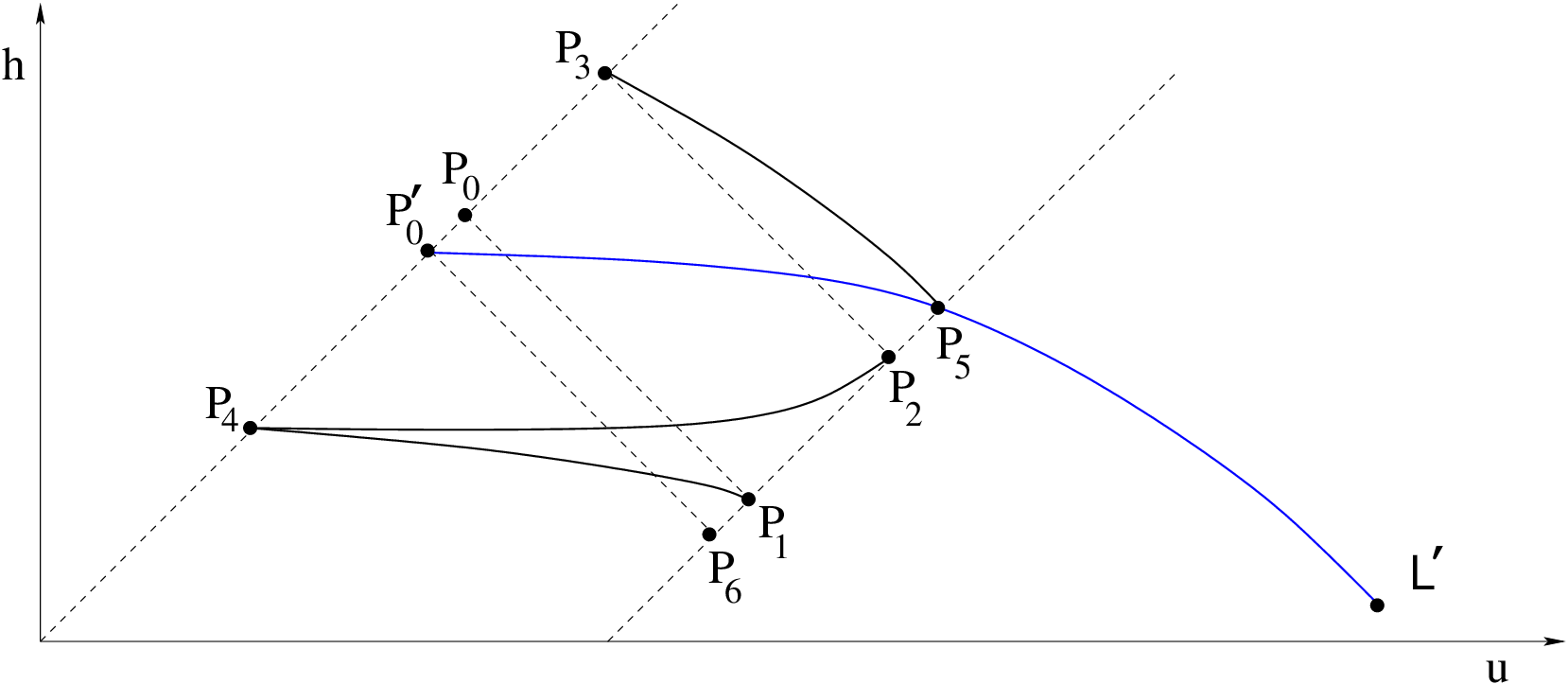}
	\caption{\small  Constructing  a perturbed 
	periodic pattern. Here the states $P_0, P_1,\ldots,P_5$ are
	the same as in Fig.~\ref{f:p96}.
	The states $P'_{0}, P_1,\ldots, P_6$ and $L'$ are those in the new pattern shown in Fig.~\ref{f:p90}.}
	\label{f:p91}.
\end{figure}

To construct our approximate solution, we begin by defining
a slightly different periodic interaction pattern.
As shown in Fig.~\ref{f:p90}, the compression wave between the states $P_1$ and $P_5$ is now split in two parts. 
\begi
\item A small portion breaks at the same point where the 1-compression
waves merge into a large 1-shock. 
This portion is thus completely cancelled by the interaction.  
\item
The remaining portion eventually forms a large 2-shock, as in the previous periodic pattern.
\endi
As a result, the new periodic approximate solution will contain an additional constant state $P_6$ between the two portions
of this 2-compression wave. 
Notice that, as $P_6\to P_1$, the new pattern becomes identical to the old one. Being able to partition the 
2-compression into two separate waves adds one more
degree of freedom to the construction of a periodic pattern.
This will be used to achieve more easily a convergence estimate.
\v
A periodic pattern as in Fig.~\ref{f:p90},
can be obtained by a perturbation argument, starting
with the pattern constructed in the previous section,
and using the implicit function theorem.
Our construction is better explained with the aid of 
 Fig.~\ref{f:p91}.

We start from our earlier periodic example in Fig.~\ref{f:p68} including states $P_0,\ldots,P_5$ and $L$.  Furthermore, by 
the  proof of Lemma 2, we have
\[h(P_0)>h(P_5).\]
We can thus find a state $P'_{0}$ along the line segment
$P_4P_0$, such that 
\bel{5.2n}
h(P'_{0})~>~h(P_5).
\eeq
Note that $P'_{0}$ can be chosen arbitrarily  close to $P_0$,

We then call $P_6$ the intersection between the 1-wave curve 
through  
$P'_{0}$ and the 2-wave curve through $P_1$. 
Notice that the 2-wave with left state $P_1$ and right state $P_6$ is a compression wave.

Finally, by (\ref{5.2n}) and using the same argument as in 
the proof of Lemma 2, we can find a new left state $L'$ which is connected to both $P'_{0}$  and $P_{5}$ by a 1-shock.

One now checks that the states $P'_0,P_1,\ldots,P_6$ and $L'$ produce the periodic pattern in Fig.~\ref{f:p91}.
In particular, notice that 
the Riemann problem with left and right states 
$P_1, P_3$ is still solved by the 1-shock $P_1 P_4$ and the 2-rarefaction $P_4 P_3$.
\v
\subsection{A sequence of nearly periodic patterns.}

\begin{figure}
	\centering
	\includegraphics[scale=.45]{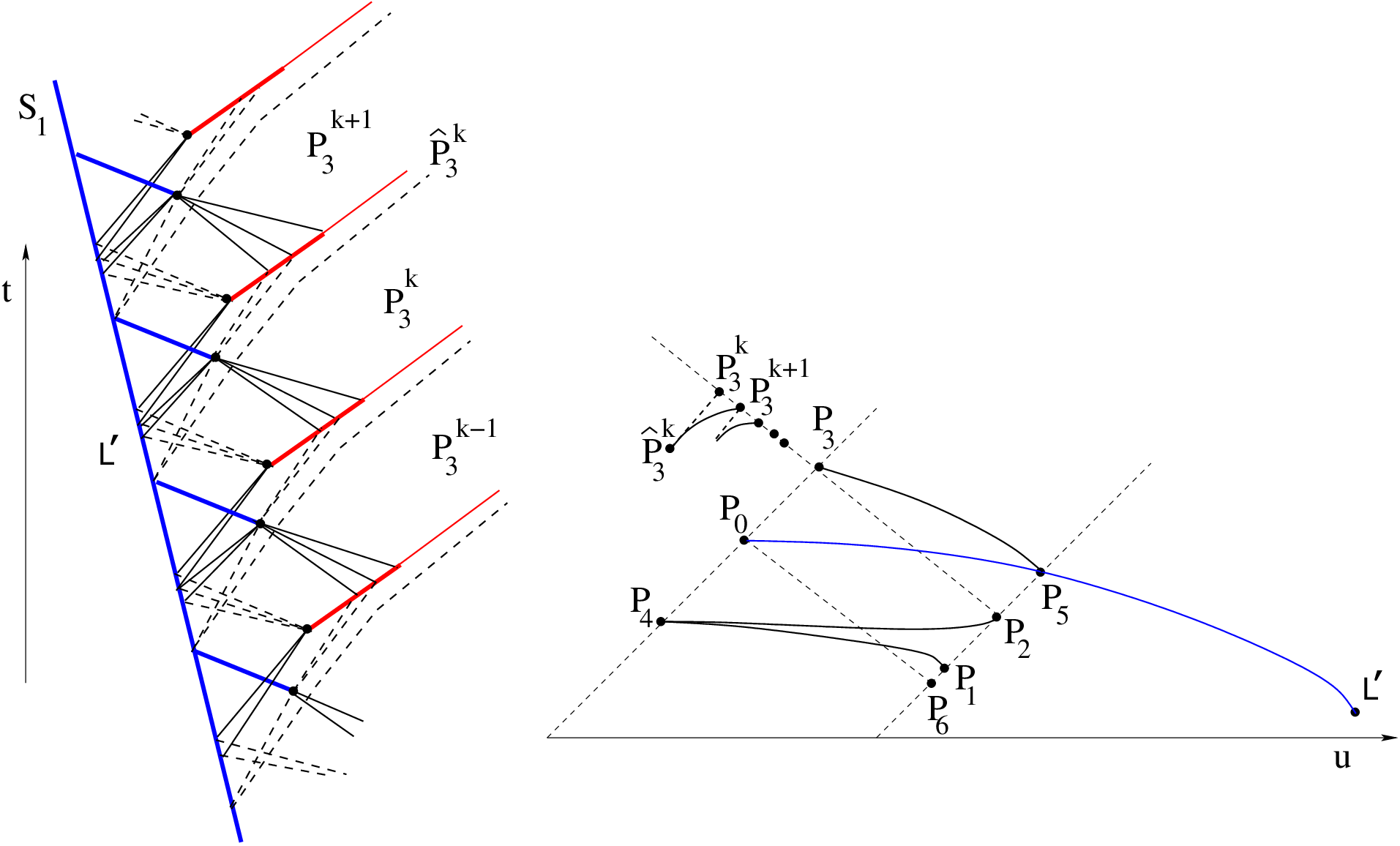}
	\caption{\small  Left: by slightly perturbing the periodic
	pattern in Fig.~\ref{f:p90}, we obtain a new 
pattern where, at the $k$-th iteration, an additional  pair of shock and rarefaction waves is produced, each with strength $\alpha/k$.   Right:  the sequence of left, middle, and  right states $P_3^{k+1}$, $\Hat P_3^k$, and $P_3^k$. }
	\label{f:p92}.
\end{figure}

Next, we slightly modify the previous periodic interaction pattern by assuming that, at the $k$-th iteration, 
the 2-shock  $P_2P_4$ is not entirely cancelled
by the 2-rarefaction $P_4 P_3$.   Instead, a pair of
2-waves survive, namely (see Fig.\ref{f:p92}):
\begi
\item a 2-shock of strength $\alpha /k$, joining the left 
state $P_3^{k+1}$ with an intermediate state $\Hat P_3^k$,
\item  a 2-rarefaction, also  of strength $\alpha /k$,
joining the  intermediate state $\Hat P_3^k$ to the right 
state~$P^k_3$.\endi
Since we require that these two shock and rarefaction waves have exactly the same strength (measured in Riemann invariants), 
all states $P^k_3$, $k=1,2,\ldots$ 
 must lie along the same 1-wave curve through $P_3$.  
Recalling that shock and rarefactions curves coincide up to second order \cite{Bbook, Sm},  for some constant $C$ we have
\bel{P3k}|P_3^{k+1}-P_3^k|~\leq~{C\over k^3}\,.\eeq
 Therefore, by choosing $\alpha>0$ small enough,
we can uniquely determine the states $P_3^k$ so that 
\bel{p3k}P_3^k~\to~P_3\qquad\qquad\hbox{as}\quad k\to\infty.\eeq
\v
In turn, we claim that  all other intermediate states $P_i^k$, with $i\in \{0,1,2,4,5,6\}$ and $k\geq 1$, can be uniquely determined as well. Indeed, these can be constructed in the following order:
\[
P_3^k~\rightarrow ~P_5^k~\rightarrow ~P_2^k~\rightarrow ~P_0^k
~ \rightarrow~ P_4^k~\rightarrow ~P_1^k~\rightarrow~ P_6^k\,.
 \]
 \begin{itemize}
 \item  $P_5^{k}$ is the state at the intersection of  the 
 1-shock curve with right state $P_3^{(k)}$  and the 1-shock curve  through $L'$.  
 \item  $P_2^k$ is the state at intersection of the 1-wave curve through $P_3$ and the 2-wave curve through $P_5^k$. 
\item  $P_0^k$ is the state at the intersection  
of the  2-wave curve through $P_3^k$ and the 1-shock curve through  $L'$.  
    \item  $P_4^k$ is the state at the intersection  of 
    the 2-shock curve with right state 
  $P_2^{k+1}$ and the 2-wave curve through $P_3^k$.      
  \item  $P_1^k$ is the state at the intersection  between the 2-wave curve through 
  $P_5^k$ and the 1-shock curve with right state $P_4^k$. 
  \item  Finally, $P_6^k$ is the state at the intersection  between the 1-wave curve through 
  $P_0^k$ and the 2-wave curve through $P_1^{k+1}$.
   \end{itemize}
We observe that, by choosing $\alpha>0$ small, 
all points $P_3^k$
will lie in a suitably small neighborhood of $P_3$.
By the implicit function theorem, all the states $P_5^k$, 
$k\geq 1$,
are well defined and lie in a suitably small neighborhood of $P_5$.

In turn, again by the implicit function theorem, it follows that
all the states $P_2^k$ are well defined and lie in a small neighborhood of $P_2$.

After six steps, all sequences of points $P^k_3 ,P^k_5,P^k_2,P^k_0,P^k_4,P^k_1,P^k_6$, $k=1,2,\ldots$
are thus uniquely determined, provided that $\alpha>0$ was chosen sufficiently small.  Moreover, we have the convergence
\bel{pikk}
\lim_{k\to\infty}P_i^k~=~P_i\qquad\qquad i=0,1,\ldots,6.\eeq
We remark that, by the convergence $P_1^k\to P_1$ and $P_6^k\to P_6$, it follows that (by possibly shrinking the value 
of $\alpha$)  the  states $P_6^k$ and $P_1^{k+1}$ are  
always connected by a 2-compression (not a 2-rarefaction).

The previous analysis achieves the construction of the modified interaction pattern shown in Fig.~\ref{f:p92}.

\subsection{An approximate solution with finite 
time BV blow-up.}
The approximate solution constructed in the previous step
(Fig.~\ref{f:p92}) contains a sequence of 2-shocks followed
by a 2-rarefaction, both of strength $\alpha/k$, $k=1,2,\ldots$.
Clearly, the total strength of all these waves is infinite.

To provide an example where blow up of the total variation
occurs in finite time, it suffices to slightly modify
the wave speeds,   so that the interaction cycles
repeat over shorter and shorter time intervals $[\tau_k,\,
\tau_{k+1}]$, with $\tau_k\to T$ as $k\to\infty$.

\begin{figure}
	\centering
	\includegraphics[scale=.5]{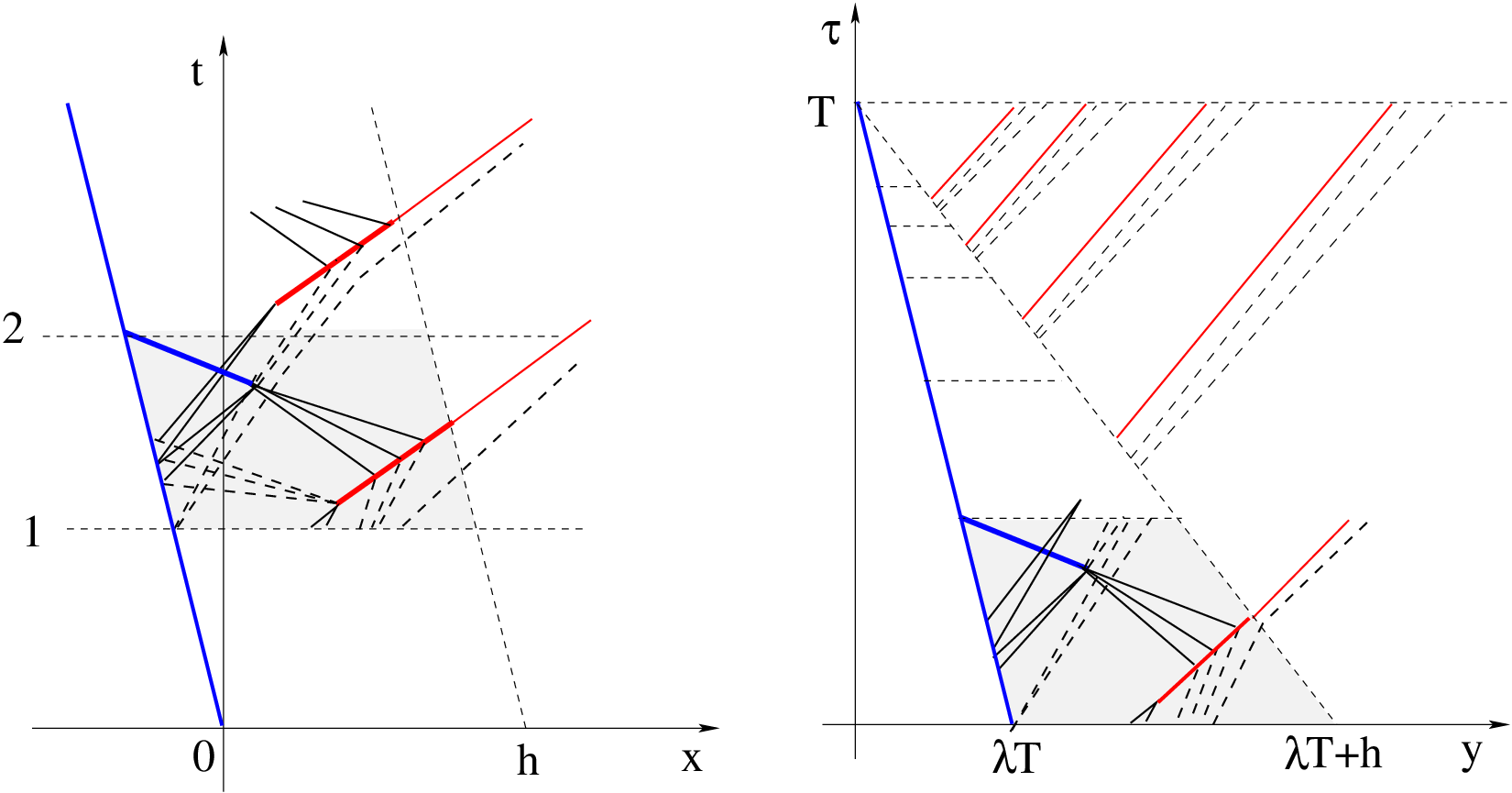}
	\caption{\small After the transformation  of the $t$-$x$ coordinates defined at (\ref{tauy}),
	from the interaction pattern shown on the left
	(same as the one in Fig.~\ref{f:p92}) we obtain a new approximate solution where all wave 
	interactions take place within  the time interval $[0,T]$.
	At the terminal time $\tau=T$, the total
	variation becomes infinite.}
	\label{f:p93}.
\end{figure}

To fix the ideas, assume that in the previous construction 
the basic interaction cycle takes place on the 
parallelograms
$$\Gamma_k~=~\Big\{(t,x)\,;~~t\in [k,\,k+1],~~x\in [-\lambda t,
-\lambda t+h]\Big\},\qquad\qquad k=1,2,\ldots$$
for some $h>0$ (see Fig.~\ref{f:p93}, left).
Fix $T>0$ sufficiently large and define
$$\ve~\doteq~-\ln\left(1-{1\over T}\right)~=~\ln\left(1+{1\over T-1}\right).$$
Observe that, as $T\to +\infty$, we have
$$\ve T~=~1+\O(1)\cdot T^{-1}.$$
Consider the transformation (see Fig.~\ref{f:p93})
\bel{tauy}\tau~=~(1-e^{-\ve t})\,T\,,\qquad\qquad y~=~(x+\lambda t+\lambda T) e^{-\ve t}\,,\eeq
defined for $t\geq 0$, $x\in\R$.  

If now $x=\xi(t)$ is the equation of a wave front in the
$t$-$x$ coordinates, let $y=\zeta(\tau)$ be the corresponding equation in the $\tau$-$y$ coordinates.
Differentiating w.r.t.~$t$ the identity
$$y(t,\xi(t))~=~\zeta\bigl(\tau(t,\xi(t))\bigr),$$
we compute
\bel{wspeed}\bega{rl}\ds\zeta'(\tau)&\ds=~{y_t + y_x \xi'\over \tau_t + \tau_x\xi'(t)}
~=~{\lambda e^{-\ve t} - (x+\lambda t + \lambda T) \ve e^{-\ve t}
+e^{-\ve t} \phi'\over T\ve e^{-\ve t}}\\[5mm]
\ds &\ds=~{\lambda(1-\ve T) + \ve (x+\lambda t) +\xi'(t)
\over \ve T}~=~\xi'(t)+\O(1)\cdot T^{-1}\,,\enda\eeq
as long as $x+\lambda t\in [0,h]$.  In other words, by 
choosing the blow up time $T$ large enough, the speeds of all 
waves contained in the strip $\{x+\lambda t\in [0,h]\}$
are almost unchanged by the coordinate transformation.
Furthermore, we impose that each pair of 2-shocks and 2-rarefactions (created at each interaction cycle) 
travels with the same speed in the $\tau$-$y$ as in the 
old $t$-$x$ coordinates (see Fig.~\ref{f:p93}), right).
 
In view of (\ref{wspeed}), this
approximate solution in the 
$\tau$-$y$ variables satisfies all conditions {\bf (C1)--(C4)}
stated at the beginning of this section.
\v

\section{A necessary condition for blowup}
\label{sec:5}
\setcounter{equation}{0}

In the second part of this paper, 
we prove that, if the total variation blows up in finite time, 
then the solution must contain an infinite number
of large shocks in a neighborhood of some point 
$\ov P=(T,\bar x)$.

Since the p-system with $p(v)= v^{-\gamma}$ 
admits a group of rescalings, a precise meaning 
of ``large shock" must be given in terms of the upper 
and lower bounds on the gas density $\rho$.
We recall that, as proved in \cite{Hoff}, for any $b>0$ 
the domain
\bel{Db}
\D_b~\doteq~\{(h,u)\,;~~h\geq 0, ~~|u|\leq b- h\}\eeq
is positively invariant for the system (\ref{1}), 
with $p$ and $h$ as in (\ref{2}), (\ref{xi_eq}).  
In the following we shall assume that the density 
$\rho$ remains uniformly positive, and hence the same 
holds for $h$. We thus consider a solution taking 
values in the domain
\bel{Dab}
\D_{ab}~\doteq~\{(h,u)\,;~~|u|\leq b- h\,,~~h\geq a\},\eeq
for some $0<a<b$.
In the following, the total variation of the vector-valued 
function $x\mapsto \bigl(h(t,x), u(t,x)\bigr)$ 
on an (possibly unbounded) interval $I\subset\R$ 
is defined as
$$\TV\bigl\{(h,u)(t,\cdot)\,;~~I \bigr\}~\doteq~
\sup\sum_i\Big(
\bigl|h(t,x_i)-h(t, x_{i-1})\bigr|+\bigl|u(t,x_i)-u(t, x_{i-1})
\bigr|\Big),$$
where the supremum is taken over all finite increasing
sequences of points $x_0<x_1<\cdots<x_N$ contained in $I$.

Observe that, as long as the solution takes values inside the compact domain $\D_{ab}$, 
a bound on the total variation of $(v,u)$
is equivalent to a bound on the total variation of 
$(h(v),u)$.   In turn, this is also equivalent 
to a bound on the total strength of waves, measured in Riemann
invariants, as in (\ref{strength}).
\v
\begin{theorem}\label{thm1}
For any two constants $b>a>0$,
there exists $\delta_0>0$ such that the following holds.
Consider an entropy weak solution $(v,u)$ 
of (\ref{1}) such that $(h(v), u)\in \D_{ab}$ for all $t,x$,
and assume that the total variation is initially bounded but
blows up at a  finite time $T$.   
Then there exists a point $\bar x$ such that
every neighborhood of $(T,\bar x)$ in the $t$-$x$ plane contains infinitely many shocks with strength $\geq \delta_0$.
\end{theorem}
\v
A proof of this theorem will be completed 
in the next two sections.
We observe that, since the initial data $(\bar v,\bar u)$
have bounded variation, for every 
$\ve_0>0$ there exists  $R_0>0$ sufficiently large such that
$$\TV\bigl\{(\bar v,\bar u)\,;~~\,]-\infty,\,-R_0[\, \bigr\}
~<~\ve_0\,.$$
$$\TV\bigl\{(\bar v,\bar u)\,;~~\,]R_0,\, +\infty[\, \bigr\}
~<~\ve_0\,,$$
For a solution
taking values in the domain $\D_{ab}$, the characteristic 
speeds $\pm c$ in  
remain uniformly bounded above and below. Indeed, since  
$a\leq h\leq b$, by  
(\ref{speed}) and (\ref{xi_eq})  
it follows
\bel{spb}c~=~\sqrt{-p'(v)}~=~\sqrt{A/\gamma}\cdot
v^{-(\gamma+1)/2}~\leq~\sqrt{A/\gamma}\cdot (b/B)^{(\gamma+1)/(\gamma-1)}
\doteq~  \hat\lambda\,.
\eeq
By choosing $\ve_0$ small enough, by the Glimm interaction
estimates it follows that for any $t>0$
the total variation of the solution
 on the two  domains
\bel{away}]-\infty\,,~-R_0-\hat \lambda t[\,,\qquad\qquad 
]R_0+\hat \lambda t\,,~+\infty[\,,
\eeq
remains uniformly small.  Here $\hat \lambda$ is the upper bound
on all characteristic speed, computed at (\ref{spb}). 
Hence, if the total variation 
blows up at time $T$,  this must happen within 
the compact interval 
$[-R_0-\hat \lambda T\,,~ R_0+\hat \lambda T]$.

We conclude this section with a preliminary lemma.
\v
\begin{lemma}\label{lemma3}
 {\it For any BV function $x\mapsto (v(x), u(x))$  with 
\bel{vb1} \bigl(h(v(x)), u(x)\bigr)\in \D_{ab}
\qquad\qquad\forall x\in\R\,,
\eeq
 the following holds.
\begi
\item[(i)] For any compact  interval $I$, one has
\bel{rw}\hbox{\rm [total strength of all waves in $I$]}~
\leq~2\Big( b+ 
\hbox{\rm [total strength of rarefaction waves in $I$]\Big).}\eeq
\item[(ii)] There exists $\delta_1>0$ such that, for every subinterval
$J\subset \R$ of length $\leq\delta_1$ one has
\bel{tw}\hbox{\rm [total strength of all waves in $J$]}~
\leq~2b\,.\eeq
\endi
}
\end{lemma}
\v
{\bf Proof.} To prove (i),
call $\mu_1,\mu_2$ the corresponding wave measures, defined as 
in Section~\ref{s:wm}. Moreover, call $\mu_i^+$  the positive part of $\mu_i$. 
Then (\ref{rw}) means that
\bel{rw2}
|\mu_1|(I)+|\mu_2|(I)~\leq~2\Big( b+ \mu_1^+(I)+\mu_2^+(I)\Big).\eeq
To prove (\ref{rw2}) we observe that the $u$ component has a downward jump at 
every point of shock.   More precisely, 
for an $i$-shock located at a point $\bar x$, one has 
$$u(\bar x +) - u(\bar x -) ~\leq~ w_i(\bar x +) - w_i(\bar x -)~<~0.$$
If $I=[\alpha,\beta]$,  since $u$ takes values inside $\D_{ab}$ we have the inequalities
$$-2b~\leq~u(\beta+)-u(\alpha-) ~\leq~\mu_1(I)+\mu_2(I),$$
$$|\mu_1| (I) -   2 \mu_1^+(I) +|\mu_2|(I) -  2  \mu_2^+(I)~\leq~2b\,.$$
This yields (\ref{rw2}).
\v
Next, if (ii) fails, then we can find a sequence of intervals 
$J_n=[\alpha_n, \beta_n]$ with lengths $\beta_n-\alpha_n=n^{-1}$,
such that 
$$|\mu_1|(J_n)+|\mu_2|(J_n)~>~2b$$
for every $n\geq 1$.   By taking a subsequence we can assume the convergence
$\alpha_n\to \bar x\in I$.   This implies 
\bel{rp3}|\mu_1|(\{\bar x\})+|\mu_2|(\{\bar x\})~\geq~2b\,.\eeq
As defined in Section~\ref{s:wm},  the left hand side of (\ref{rp3}) 
is the total strength of the two waves in the solution of the 
Riemann problem with left and right states
$(v,u)(\bar x-)$, $(v,u)(\bar x+)$.
Since this solution takes values in the domain $\D_{ab}$, the sum of these
two strengths must be $\leq 2(b-a)$.  We thus reach a contradiction with (\ref{rp3}),
proving the second part of the lemma. 
\endproof
\v
\section{Wave  decay estimates} 
\label{sec:7}
In this section we prove two estimates on the decay of rarefaction waves, extending the analysis in Section~3
to general BV solutions.

\begin{figure}[htbp]
\centering
  \includegraphics[scale=0.4]{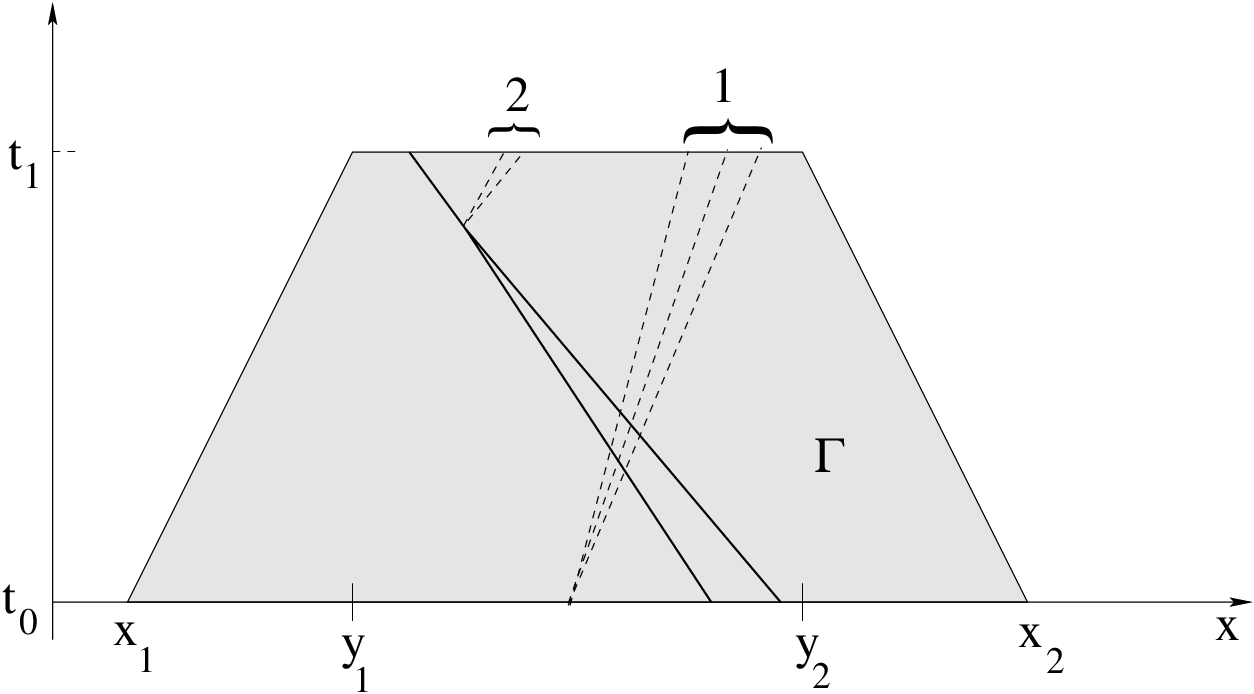}
    \caption{\small  The two types of rarefaction waves 
    which can cross the upper boundary
     of the trapezoidal domain $\Gamma$ at (\ref{Gd1}).   }
\label{f:hyp122}
\end{figure}
\subsection{Solutions without large shocks.}
We first study the simpler case where no large shock is present.
Consider a BV solution of (\ref{1}).  Fix a time step $\Delta t$ and a space step 
$\Delta x$, such that
\bel{Dxt}
\Delta x~=~2\hat\lambda \cdot \Delta t\,,\eeq
and consider a
domain of the form (see Fig.~\ref{f:hyp122})
\bel{Gd1}
\Gamma~=~\Big\{ (t,x)\,;~~t\in [t_0, t_1],\quad
x\in [ x_1+\hat\lambda
(t-t_0)\,,~ x_2-\hat\lambda
(t-t_0)]\Big\},\eeq
with  
\bel{t01} t_1-t_0~=~\Delta t\,,\qquad\qquad x_2-x_1~=~2\Delta x\,.
\eeq
We seek an estimate on the total amount of rarefaction waves
at time $t_1$, along the upper boundary 
$$[y_1,y_2]~=~ \bigl[ x_1+\hat\lambda
(t_1-t_0)\,,~ x_2-\hat\lambda
(t_1-t_0)\bigr].$$
Notice that, by (\ref{Dxt}) and (\ref{t01}), this interval has length 
$y_2-y_1=\Delta x$.

 As shown in Fig.~\ref{f:hyp122}, these rarefactions 
 can be of two types:
\begi
\item[1)] Old rarefactions which were already present along the 
basis $[x_1,x_2]$ at time $t_0$.   The total amount of these
waves can be controlled because, as discussed in Section~3,
their density has decayed during the entire time interval 
$[t_0, t_1]$. 
Roughly speaking, we have
$$\bega{l}
\hbox{[total amount of old rarefactions]}~\leq~(y_2-y_1)\cdot
\hbox{[maximum density]}\\[4mm]
\qquad \ds \leq~(y_2-y_1)\cdot
{\O(1)\over t_1-t_0}~=~\O(1)\cdot 2\hat\lambda\,.\enda$$

\item[2)] New rarefactions produced by wave interactions inside the domain $\Gamma$.  Assuming that the total strength of all waves at the initial time $t_0$ is $\leq K_0$ and all shocks in 
$\Gamma$ have size $\leq\delta_0$, the total strength of 
these new waves will be of order $\O(1)\cdot \delta_0 K_0^2$.
\endi
\v
\begin{lemma}\label{lemma4}
{\it  Let $\D_{ab}$ be the domain
in (\ref{Dab}).  Then one can find a  constant $K_{ab}$ 
such that,
for any given  $K_0$, there exists 
$\delta_0>0$ for which the following holds. 

Let $(v,u)$ be a BV solution of (\ref{1}) taking values inside 
$\D_{ab}$ and let $\Gamma$ be the trapezoid defined at (\ref{Gd1}).  
Assume that
\begi
\item[(i)] at time $t_0$ the total strength of all waves
contained inside the lower boundary $[x_1, x_2]$ is $\leq K_0$, and
\item[(ii)] all shocks inside $\Gamma$ have strength 
$\leq\delta_0$.
\endi
Call $\hat \mu_i$ the measures of $i$-waves in the solutions at time $t_1$.
Then the total strength of all rarefactions waves
contained inside  the upper boundary $[y_1,y_2]$ satisfies 
\bel{rar1} (\hat \mu_1^++\hat \mu_2^+)\bigl([y_1,y_2]\bigr)~\leq~ K_{ab}\,.\eeq
}
\end{lemma}
Notice that here the constant $K_{ab}$ can be large, but is independent of 
$K_0$.   This implies that, if the initial data contain a large amount of waves but the 
solution does not develop large shocks, then most of the rarefaction
waves present at time $t_0$ will disappear during the time interval $[t_0,t_1]$, being canceled by 
waves of the same family but opposite sign. Compared with the decay estimate proved in Chapter 10 of 
\cite{Bbook},
the main difference is that here the total 
strength of waves can be large. However, 
thanks to Proposition~\ref{prop1}, 
the total strength of new waves
produced by interactions can be made arbitrarily small  by 
choosing $\delta_0>0$ small enough.
\v
{\bf Proof.}  {\bf 1.} Consider any interval $I\doteq [a, b]~\subset~[y_1,y_2]$.
Call $t\mapsto \xi_1(t)$, $t\mapsto \xi_2(t)$ respectively the
minimal backward
2-characteristics passing through $a, b$ at time $t_1$.
Motivated by (\ref{ch3}), we define
\bel{Zdef4}Z(t)~\doteq ~\int_{\xi_1(t)}^{\xi_2(t)}
\vp(t,x)\, dx\,,\eeq
where 
 $$\vp(t,x)~=~\vp\bigl(w_2(t,x)-w_1(t,x)\bigr)
~=~\vp(2h(t,x))$$  
is the function defined at (\ref{eqphi}).
As long as the solution $(v,u)$ takes values in $\D_{ab}$ we have 
$h(v)\geq a$, hence
$\vp$ remains bounded and uniformly positive.
The  integral (\ref{Zdef4}) thus 
provides an equivalent way to 
measure the distance
$\xi_2(t)-\xi_1(t)$ between the two characteristics.

Since the total variation is bounded and all characteristic speeds are bounded by 
$\hat\lambda$, the function
$t\mapsto Z(t)$ is Lipschitz continuous. Its time derivative
\bel{DZ}
\dot Z~=~\dot\xi_2 \cdot\vp(\xi_2) - \dot\xi_1 \cdot\vp(\xi_1) 
+\int_{\xi_1}^{\xi_2}D_t\vp\eeq
is well defined for a.e.~time $t$. Notice that here 
$D_t\vp$ is a bounded measure. Its atomic part is supported on the set of shocks. 


In the following, for any given time $t$ we denote by $\S$  the set of all shocks contained inside
the interval $[a(t), b(t)]$ and call
$k_\alpha\in \{1,2\}$ the family of the  shock 
located at $x_\alpha(t)$. Moreover, by $D^c_x$ be denote the continuous (i.e., non atomic) part of a distributional derivative w.r.t.~$x$, 
Motivated by (\ref{ch3}), denoting by $c(x)=c(v(t,x))$ the wave speed and using 
(\ref{w12}), (\ref{DZ}), we compute
\bel{DZ2}\bega{rl}
\dot Z&\ds=~\ds c(\xi_2)\vp(\xi_2) -c(\xi_1)
\vp(\xi_1)  - \int_{\xi_1}^{\xi_2}\vp'\, c\cdot 
\bigl( D_x^c w_1 +D_x^c w_2)
 \\[4mm]&\qquad\ds
 -  \sum_{\alpha\in \S} 
\dot x_\alpha\cdot
\bigl[ \vp(x_\alpha+)- \vp(x_\alpha-)\bigr]\\[4mm]
\ds
 \qquad&\ds=~ \int_{\xi_1}^{\xi_2}(c'\vp+c\vp')\cdot 
\bigl( D_x^c w_2 -D_x^c w_1) + \sum_{\alpha\in\S}
\bigl[ c(x_\alpha+)\vp(x_\alpha+)-c(x_\alpha-) \vp(x_\alpha-)\bigr]
 \\[4mm]
 &\qquad \ds - \int_{\xi_1}^{\xi_2}\vp'\,c\cdot 
\bigl( D_x^c w_1 +D_x^c w_2) -
   \sum_{\alpha\in \S} 
\dot x_\alpha\cdot
\bigl[ \vp(x_\alpha+)- \vp(x_\alpha-)\bigr]. \enda
\eeq
For each shock $\alpha\in\S$, two cases must be considered.

CASE 1: The shock at $x_\alpha$ belongs to the first family.
By definition, its strength is 
$$\sigma_\alpha~=~w_1(x_\alpha-)-w_1(x_\alpha+)~\in~[0,\delta_0].$$
In this case, we have
$$\dot x_\alpha~=~-{c(x_\alpha+)+c(x_\alpha-)\over 2 }+ \O(1)\cdot \sigma_\alpha^2\,,$$
$$\bigl|w_2(x_\alpha+)-w_2(x_\alpha-)\bigr|~=~\O(1)\cdot\sigma_\alpha^3
\,,$$
$$\vp(x_\alpha+)-\vp(x_\alpha-)~=~\vp'(x_\alpha-)\sigma_\alpha +\O(1)\cdot
\sigma_\alpha^2\,,$$
$$c(x_\alpha+)-c(x_\alpha-)~=~c'(x_\alpha-)\sigma_\alpha +\O(1)\cdot
\sigma_\alpha^2\,,$$
Using  the fundamental relation
$c'\vp+2c\vp'~=~0$ we thus obtain
\bel{sh1}\bega{l}
\bigl[ c(x_\alpha+)\vp(x_\alpha+)-c(x_\alpha-) \vp(x_\alpha-)\bigr]-
\dot x_\alpha\cdot
\bigl[ \vp(x_\alpha+)- \vp(x_\alpha-)\bigr]\\[4mm]
\qquad =~\bigl[ c'(x_\alpha-) \vp(x_\alpha-) + c(x_\alpha-) \vp'(x_\alpha-) \bigr] 
\,\sigma_\alpha + c(x_\alpha-)\vp'(x_\alpha-) \sigma_\alpha +\O(1)\cdot \sigma_\alpha^2\\[4mm]
\qquad =~\O(1)\cdot \sigma_\alpha^2\,.\enda
\eeq

CASE 2: The shock at $x_\alpha$ belongs to the second family.
By definition, its strength is 
$$\sigma_\alpha~=~w_2(x_\alpha-)-w_2(x_\alpha+)~\in~[0,\delta_0].$$
In this case, we have
$$\dot x_\alpha~=~{c(x_\alpha+)+c(x_\alpha-)\over 2 }+ 
\O(1)\cdot \sigma_\alpha^2\,,$$
$$\bigl|w_1(x_\alpha+)-w_1(x_\alpha-)\bigr|~=~\O(1)\cdot\sigma_\alpha^3
\,,$$
$$\vp(x_\alpha+)-\vp(x_\alpha-)~=~-\vp'(x_\alpha-)\sigma_\alpha +\O(1)\cdot
\sigma_\alpha^2\,,$$
$$c(x_\alpha+)-c(x_\alpha-)~=~-c'(x_\alpha-)\sigma_\alpha +\O(1)\cdot
\sigma_\alpha^2\,,$$
In this case we  obtain
\bel{sh2}\bega{l}
\bigl[ c(x_\alpha+)\vp(x_\alpha+)-c(x_\alpha-) \vp(x_\alpha-)\bigr]-
\dot x_\alpha\cdot
\bigl[ \vp(x_\alpha+)- \vp(x_\alpha-)\bigr]\\[4mm]
\qquad =~\bigl[ c'(x_\alpha-) \vp(x_\alpha-) + c(x_\alpha-) \vp'(x_\alpha-) \bigr] 
\,\sigma_\alpha - c(x_\alpha-)\vp'(x_\alpha-) \sigma_\alpha +\O(1)\cdot \sigma_\alpha^2\\[4mm]
\qquad =~ c'(x_\alpha-) \vp(x_\alpha-)  \sigma_\alpha+\O(1)\cdot \sigma_\alpha^2\,.\enda
\eeq
{}From   (\ref{DZ2}), 
using (\ref{sh1})-(\ref{sh2}) and  the relation 
$c'\vp+2c\vp'=0$ one obtains
\bel{DZ3}
\dot Z~=~
 \int_{\xi_1}^{\xi_2}c'\vp\cdot  D_x^c w_2  + \sum_{\alpha\in\S_2}
 c'(x_\alpha-)\vp(x_\alpha-)\,\sigma_\alpha +\O(1)\cdot \sum_{\alpha\in\S}
\sigma_\alpha^2\,.
\eeq
Here the first summation ranges over the set $\S_2$ of all shocks of the 
second family, while the second summation ranges over the set of all shocks (of both families).
\v
{\bf 2.} 
Call $\hat\mu_2$ the measure of 2-waves
in the solution at time $t_1$. For any $\ve>0$
we can find finitely many intervals
$[a_\ell, b_\ell]$, $\ell=1,\ldots,m$, whose union contains
nearly all positive 2-waves, and very few negative 2-waves.
More precisely:
\bel{mu3}
\hat\mu_2^+\Big([y_1,y_2]\setminus 
\bigcup_{\ell}[a_\ell, b_\ell]
\Big)~\leq~\ve,\qquad\qquad 
\hat\mu_2^-\Big( \bigcup_{\ell}[a_\ell, b_\ell]\Big)~\leq~\ve.\eeq
For each such interval, let $\xi_\ell(t), \tilde\xi_\ell(t)$ be the 
minimal backward 2-characteristics through
$a_\ell, b_\ell$, respectively.  
Setting
$$Z_\ell(t)~\doteq~\int_{\xi_\ell(t)}^{\tilde\xi_\ell(t)}
\vp(t,x) \, dx$$
and applying (\ref{DZ3}) to each subinterval $[\xi_\ell, \tilde\xi_\ell]$
we obtain
\bel{DZ4}
\sum_\ell \dot Z_\ell(t)
~=~
 \int_{\xi_\ell}^{\tilde \xi_\ell }c'\vp\cdot  D_x^c w_2  + \sum_{\alpha\in\S_2}
 c'(x_\alpha-)\vp(x_\alpha-)\,\sigma_\alpha + \O(1)\cdot\sum_{\alpha\in\S}
\sigma_\alpha^2\,,\eeq
where now $S_2$ and $\S$ refer to the shocks contained in the union of the intervals $[\xi_\ell(t),\tilde \xi_\ell(t)]$. 
For convenience, we  introduce the constants
$$0~<~\kappa_{min}~\doteq~\min\, c'\vp\,,\qquad\qquad
 \kappa_{max}~\doteq~\max \,c'\vp\,,\qquad\qquad\vp_{max}~\doteq~\max\,\vp\,,$$
defined by taking the minimum and the maximum values of the functions 
$c'\vp$ and $\vp$ over the domain $\D_{ab}$.  
Using Lemmas~\ref{lemma5} and \ref{lemma7}, 
the amounts of positive and negative 2-waves
contained in the union of the intervals $[\xi_\ell(t),\tilde\xi_\ell(t)]$
at any time $t\in [t_0, t_1]$ can be estimated as
\bel{mu2+}\mu_2^+\left(\bigcup_\ell[\xi_\ell(t), \tilde\xi_\ell(t)]\right)~\geq~
\mu_2^+\left(\bigcup_\ell[\xi_\ell(t_1), \tilde\xi_\ell(t_1)]\right)-\O(1)\cdot
\delta_0\,
\mu^{int}(\Gamma),\eeq
\bel{mu2-}\mu_2^-\left(\bigcup_\ell[\xi_\ell(t), \tilde\xi_\ell(t)]\right)~\leq~
\mu_2^-\left(\bigcup_\ell[\xi_\ell(t_1), \tilde\xi_\ell(t_1)]\right)+\O(1)\cdot
\delta_0\,\mu^{int}(\Gamma).\eeq
Combining (\ref{DZ4}) with (\ref{mu2+})-(\ref{mu2-}) we obtain
\bel{DZ5}\bega{rl}\ds\sum_\ell \dot Z_\ell&\ds\geq~\kappa_{min}\cdot\mu_2^+\left(\bigcup_\ell\,[\xi_\ell(t), \tilde\xi_\ell(t)]\right) - \kappa_{max}\cdot\mu_2^-\left(\bigcup_\ell\, [\xi_\ell(t), \tilde\xi_\ell(t)]\right)  +\O(1)\cdot\sum_{\alpha\in \S} \sigma_\alpha^2\\[4mm]
&\ds\geq~\kappa_{min}\cdot \hat \mu_2^+\Big(\bigcup_\ell\,[a_\ell,\, b_\ell]\Big) - \kappa_{max}\cdot\hat \mu_2^-\Big(\bigcup_\ell\,[a_\ell,\, b_\ell]\Big)
-\O(1)\cdot \delta_0\,\mu^{int}(\Gamma)- 
\O(1)\cdot\delta_0\, V(t)\\[4mm]
&\ds\geq~\kappa_{min}\cdot\Big(\hat \mu_2^+([y_1,y_2])-\ve\Big)-
\kappa_{max}\cdot\ve
-\O(1)\cdot \delta_0\,\mu^{int}(\Gamma)- 
\O(1)\cdot\delta_0\,  K_0\\[4mm]
&\ds\geq~\kappa_{min}\cdot\hat \mu_2^+([y_1,y_2])-\O(1)\cdot \delta_0 K_0^2 - 
\O(1)\cdot\ve.\enda
\eeq
Observing that
$$ \sum_\ell Z_\ell(t_0)~\geq~ 0\,,\qquad\qquad 
\sum_\ell Z_\ell(t_1)~\leq~(y_2-y_1)\, \vp_{max}\,,$$
and integrating (\ref{DZ5}) over the time interval $[t_0, t_1]$, 
we obtain
\bel{Zm1}
(t_1-t_0)\kappa_{min}\cdot\hat \mu_2^+([y_1,y_2])~\leq~(y_2-y_1)\, \vp_{max}+\O(1)\cdot (t_1-t_0) \delta_0 K_0^2 +
\O(1)\cdot (t_1-t_0)\ve \,.\eeq
\v
{\bf 3.} 
Since $y_2-y_1= \hat \lambda (t_1-t_0)$ and $\ve>0$ can be taken arbitrarily small,
(\ref{Zm1}) yields an a priori bound
on the total amount of positive 2-waves at the terminal  time $t_1$, namely
\bel{Zm2}
\hat \mu_2^+([y_1,y_2])~\leq~\hat\lambda\,{\vp_{max}\over \kappa_{min}}
+C_1\delta_0 K_0^2\,,\eeq
for a suitable constant $C_1$.  Of course, an entirely similar estimate
is valid for rarefaction waves of the first family.

For any given $K_0$, we can now choose $\delta_0>0$ so that 
\bel{Zm3}\hat \mu_1^+([y_1,y_2])+
\hat \mu_2^+([y_1,y_2])~\leq~2\hat\lambda\,{\vp_{max}\over \kappa_{min}}
+2C_1\delta_0 K_0^2~\leq~3\hat\lambda\,{\vp_{max}\over \kappa_{min}}~\doteq~K_{ab}\,.\eeq
With the above definition of the constant $K_{ab}$,  the conclusion of the Lemma is achieved. 
\endproof
\v
\v
\subsection{Solutions with one large shock.}

Our next goal is to  extend Lemma~\ref{lemma4} to the  
case where one large shock is present.  

To fix the ideas, let $x=\gamma(t)$ be the location
of a 1-shock, in a solution $U=(v,u)$ of (\ref{1}).
For a BV solution, the local behavior near the shock is well 
understood.  The shock speed 
$t\mapsto \dot \gamma(t)$ is a BV function
with at most countably many jumps.  These occur on a countable set 
$\T$ of times where another shock (or a centered compression) impinges on 
$\gamma$.
The left and right limits of the solution across the shock
\bel{Ulr}
U^l(t)~\doteq~\lim_{x\to \gamma(t)-} U(t,x)\,,\qquad\qquad
U^r(t)~\doteq~\lim_{x\to \gamma(t)+} U(t,x)\,,\eeq
are well defined for all times $t\notin\T$. 
In terms of these limits one can define the measures
$\mu_i^l$, $\mu_i^r$ of waves to the left and to the 
right of the shock as follows.   

Call $(w_1^l(t), w_2^l(t))$ the Riemann coordinates of the 
state $U^l(t)$.   Moreover, denote by $\mu_i^{(t)}$ the measure
of $i$-waves in the solution $u(t,\cdot)$.  Then
\begi
\item The continuous (i.e., non atomic) part of $\mu_i^l$
coincides with the continuous part of $- D_t w_i^l$.
\item If $\tau\in \T$ is a time where $U^l$
has a jump, then 
$$\mu_1^l(\{\tau\})~\doteq~\lim_{\ve\to 0+}
\mu_1^{(\tau-\ve)}\Big( [\gamma(\tau),\, \gamma(\tau-\ve)[ \,\Big),$$ 
 is the amount of  1-waves hitting the shock from the left, at time $\tau$.
 Moreover, 
 


$$\mu_2^l(\{\tau\})~\doteq~\lim_{\ve\to 0+}
\mu_2^{(\tau-\ve)}\Big( [\gamma(\tau)-\hat\lambda\ve,\, \gamma(\tau)]
\Big).$$ 
is the amount of 2-waves hitting the shock from the left, at time $\tau$.

\endi
Similarly, let $(w_1^r(t), w_2^r(t))$ be the Riemann coordinates of the 
state $U^r(t)$.   Then
\begi
\item The continuous  part of $\mu_1^r$
coincides with the continuous part of $ D_t w_1^r$. 
The continuous  part of $\mu_2^r$
coincides with the continuous part of $- D_t w_2^r$.
\item If $\tau\in \T$ is a time where $U^r$
has a jump,
then 
$$\mu_1^r(\{\tau\})~\doteq~\lim_{\ve\to 0+}
\mu_1^{(\tau-\ve)}\Big( \,]\gamma(\tau-\ve),\, \gamma(\tau)+\hat\lambda\ve]
\Big).$$ 
is the amount of 1-waves hitting the shock from the right, at time $\tau$, while 
$$\mu_2^r(\{\tau\})~\doteq~\lim_{\ve\to 0+}
\mu_2^{(\tau+\ve)}\Big( [\gamma(\tau),\, \gamma(\tau)+\hat \lambda\ve]\Big),$$ 
is the amount of 2-waves coming out from the shock from the right. 
\endi
For an analysis of the local structure of a solution 
in the neighborhood of a point, we refer to
\cite{BLF, DP}
\begin{figure}[htbp]
\centering
 \includegraphics[scale=0.5]{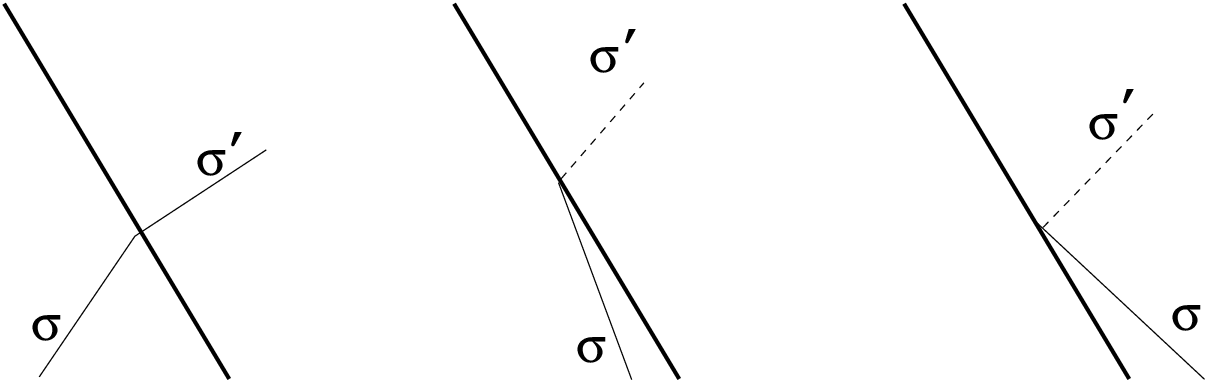}
    \caption{\small  As long as all states remain inside the domain $\D_{ab}$, the strength  of the 2-wave
    emerging from the interaction of a small wave 
    with a large 1-shock satisfies the bounds  
 in Lemma~ \ref{lemma0}.}
\label{f:p98}
\end{figure}

\begin{figure}[htbp]
\centering
 \includegraphics[scale=0.4]{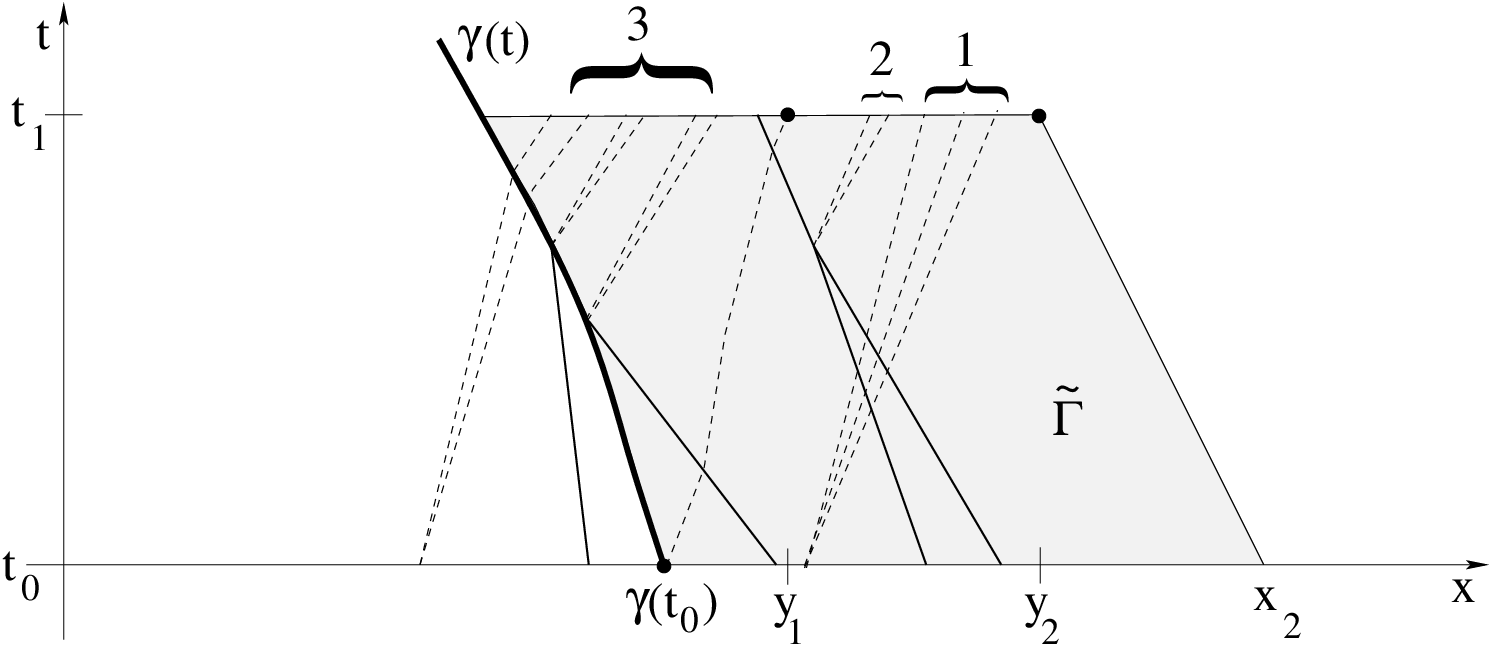}
    \caption{\small  
    The three types of rarefaction waves 
    which can cross the upper boundary
     of the large trapezoidal domain $\Gamma$ at (\ref{Gd3}).  
     Here $y_1$ is the largest point such that the minimal backward 
     2-characteristic through $(t_1. y_1)$ meets the 1-shock $\gamma$ .}
\label{f:p99}
\end{figure}

Recalling (\ref{Dxt}), we now consider the trapezoidal domain 
$\Tilde\Gamma$ shown in 
Fig.~\ref{f:p99} where the
left side is a large 1-shock.  More precisely:
\bel{Gd3}
\Tilde \Gamma~=~\Big\{ (t,x)\,;~~t\in [t_0, t_1],\quad
x\in \,]\gamma(t)\,,~ x_2-\hat\lambda
(t-t_0)]\Big\},\eeq
where now 
\bel{gsize}
t_1-t_0~=~\Delta t,\qquad\qquad x_2-\hat\lambda(t_1-t_0) -\gamma(t_1)~\leq~ 2\Delta x\,.\eeq

In order to estimate the strength of 2-waves emerging from 
interactions with the large 1-shock, we rely on the following 
elementary estimate.  

By taking limits of front tracking approximations and 
using the lower semicontinuity 
of wave measures w.r.t.~$\L^1$ convergence (proved in \cite{Bbook}),
from Lemma~\ref{lemma0} one obtains

\begin{lemma}\label{lems2}
Consider a  1-shock, located along the curve  
$\{x= \gamma(t)\,;~t\in [t_0, t_1]\}$. 
For $i=1,2$, call $\mu_i^l,\mu_i^r$ the measures of $i$-waves to 
the left and to the right
of the shock, defined as above. Then the total strength of 2-waves
emerging from the shock is estimated by
\bel{2w}
|\mu_2^r|\bigl([t_0, t_1]\bigr)~\leq~
C_\gamma\cdot\bigl(|\mu_1^l|+|\mu_2^l|\bigr)\bigl([t_0, t_1]\bigr)
+|\mu_1^r|\bigl([t_0, t_1]\bigr).\eeq
\end{lemma}

Because of the above estimates,
it is convenient to measure the  {\bf weighted strength}
of all small waves contained in $\Gamma$ at time $t$ 
by assigning a larger weight to waves
which are approaching the large 1-shock. Namely:
\bel{Wt}\bega{l}
\ds W(t)~\doteq~\bigl(2 |\mu_1|
 + |\mu_2|\bigr)\Big(] \gamma(t)\,,~ x_2- \hat\lambda(t-t_0)]\Big).\enda\eeq

As before, we seek an estimate on the total strength
of rarefaction waves at time $t_1$. 
As shown in Fig.\ref{f:p99}, three different types of rarefactions can now occur.
\begi
\item[1)] Old rarefactions, which were already present along the 
bottom side  $[\gamma(t),\,x_2]$ at time $t_0$.
Because of genuine nonlinearity, these waves decay. Their total strength at time $t_1$ is uniformly bounded, regardless of the amount of waves at the initial time $t_0$.
\item[2)] New rarefactions, produced by 
interactions of small waves inside the domain $\Tilde \Gamma$. By
the interaction estimates (\ref{299})--(\ref{22}), the total strength 
of these waves can be rendered arbitrarily small, 
by choosing $\delta_0$ small enough. 
\item[3)]  2-rarefactions emerging from the large shock.
These  occur  when (i) a 2-rarefaction crosses the large shock, or (ii)
a 1-shock (or a 1-compression wave) hits the large 1-shock 
from the left or from the right.  Because of the definition (\ref{Wt}),
the total strength of these waves
is controlled by the decay in the functional $W(\cdot)$.
\endi

\begin{lemma}\label{lemma8}
 {\it Consider the domain $\D_{ab}$ 
in (\ref{Dab}).  Then one can find a  constant $\Tilde K_{ab}$ 
such that,
for any given  $\Tilde K_0$, there exists 
$\delta_0>0$ for which the following holds. 

Let $(v,u)$ be a BV solution of (\ref{1}) taking values inside 
$\D_{ab}$ and let $\Tilde \Gamma$ be the domain in 
(\ref{Gd3})-(\ref{gsize}), 
shown in Fig.~\ref{f:p99}. Assume that:
\begi
\item[(i)] At time $t_0$ the total weighted strength of all waves
along the lower boundary $[\gamma(t_0), x_2]$ is $W(t_0)\leq \Tilde  K_0$.
\item[(ii)] 
The total strength of all waves impinging on the large 1-shock from the left
is $\leq 3b\,$. 
\item[(iii)] All shocks inside $\Tilde \Gamma$ have strength $\leq\delta_0$.
\endi
Then, calling $\hat\mu_1,\hat \mu_2$ the wave measures at time $t_1$, the total strength of all rarefactions
contained inside the upper boundary 
of $\Tilde\Gamma$ satisfies
\bel{rar4} (\hat \mu_1^++\hat \mu_2^+)\bigl([\gamma(t_1), x_2-\hat\lambda(t_1-t_0)]\bigr)~\leq~ \Tilde K_{ab}+
\bigl[W(t_0)- W(t_1)\bigr].\eeq
}
\end{lemma}
\v
{\bf Proof.}   
{\bf 1.}
As shown in Fig.~\ref{f:p99}, consider the points
$$y_1~\doteq~\inf\Big\{ y>\gamma(t_1)\,;~~\hbox{the minimal 2-characteristic
$\xi$ through $(t_1,y)$ satisfies $\xi(t_0)\geq \gamma(t_0)$}\Big\}.$$
$$y_2~\doteq~x_2-\hat\lambda(t_1-t_0)\,.$$
{}From (\ref{gsize}) we deduce
\bel{y123}
y_2-\gamma(t_1)~\leq~2\Delta x\,,\qquad\qquad x_2-\gamma(t_0)~\leq~
3\Delta x\,.\eeq

\begin{figure}[htbp]
\centering
 \includegraphics[scale=0.35]{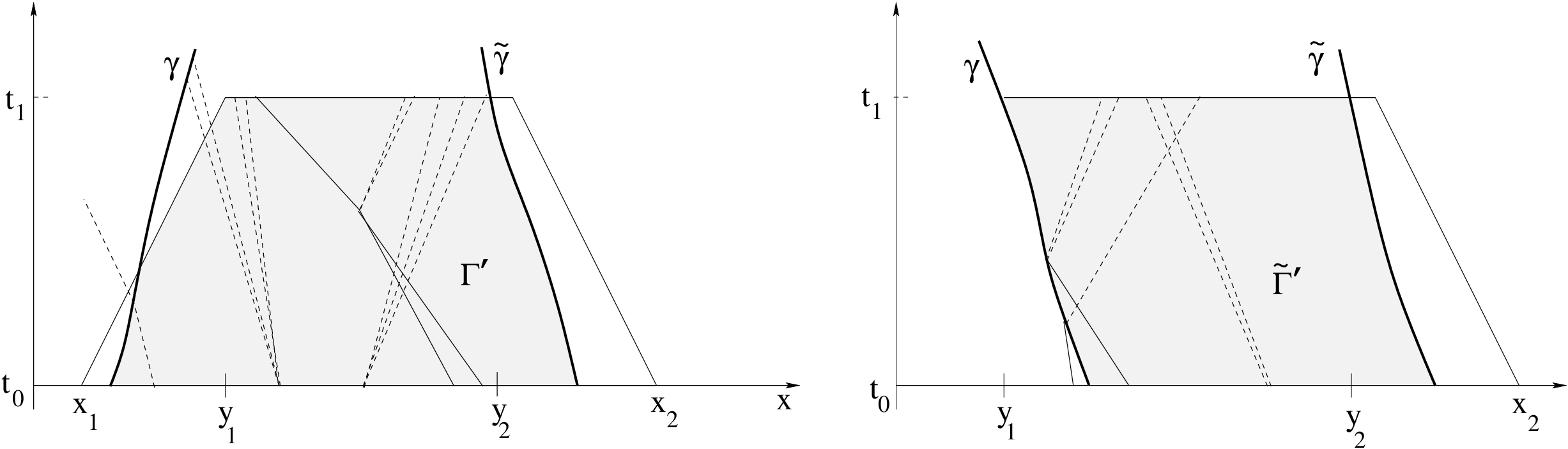}
    \caption{\small  The domains $\Gamma'$  and  $\Tilde \Gamma'$ 
    considered in Remark~\ref{remark2}. }
\label{f:p100}
\end{figure}

We observe that, at time $t_1$, all 1-rarefactions located inside the
interval $[\gamma(t_1), y_2]$ are of Type 1 (old rarefactions already 
present at time $t_0$ or Type 2 
(new rarefactions generated by wave interactions inside $\Tilde\Gamma$).
These can be estimated by the same techniques as in Lemma~\ref{lemma4}.
Similarly, at time $t_1$ all 2-rarefactions located inside the
interval $[y_1, y_2]$ are of Type 1 or 2, and can be bounded in the same way.
Calling $\hat\mu_i$ the measure of $i$-waves at time $t_1$, 
by choosing $\delta_0>0$ sufficiently small we can thus conclude
\bel{est5} 
\hat\mu_1^+\Big([\gamma(t_1), y_2]\Big)~\leq~K_{ab}\,,
\qquad\qquad \hat\mu_2^+\Big([y_1,\, y_2]\Big)~\leq~K_{ab}\,.\eeq
for a suitable constant $K_{ab}$.
\v
{\bf 2.} To prove the lemma, it thus remains to bound the total amount 
of 2-rarefactions contained in the  interval $[\gamma(t_1), y_1[\,$.
This can be achieved by standard interaction estimates.
Indeed, using (\ref{2w}) we obtain
\bel{we1}
\bega{l}
|\hat\mu_2| \Big( [\gamma(t_1), y_2[\Big)~ \leq ~[\hbox{total strength of 2-waves emerging from the big 1-shock}]\\[4mm]
\qquad\qquad\qquad +[\hbox{total strength of new waves
generated by interactions within $\Tilde \Gamma$}]\\[4mm]
\qquad\leq~|\mu_2^r|\bigl([t_0, t_1]\bigr) + \O(1)\cdot\delta_0\\[4mm]
\qquad\leq~C_\gamma\cdot\bigl(|\mu_1^l|+|\mu_2^l|\bigr)\bigl(]t_0, t_1]\bigr)
+|\mu_1^r|\bigl(]t_0, t_1]\bigr)+\O(1)\cdot\delta_0\,.
\enda\eeq
The assumption (ii) implies
\bel{we3}
C_\gamma\cdot\bigl(|\mu_1^l|+|\mu_2^l|\bigr)\bigl([t_0, t_1]\bigr)~\leq~3b\,
 C_\gamma\,.\eeq
By Lemma~\ref{lemma0}, when a 
1-wave impinges on the large 1-shock from the right,
the outgoing 2-wave resulting from the interaction has smaller strength. 
By our definition of $W$, the incoming 1-wave has weight 
2 while
the outgoing 2-wave has weight 1.  Keeping this in mind, we obtain
\bel{Wt01}\bega{l}
W(t_1)-W(t_0)~\leq~C_\gamma\cdot [\hbox{total strength of waves impinging on the big 1-shock
from the left}]\\[4mm]
\qquad -[\hbox{total strength of 1-waves impinging on the big 1-shock from the right}] \\[4mm]
\qquad+ \O(1)\cdot[\hbox{total strength of new waves
generated by interactions within $\Tilde \Gamma$}]\\[4mm]
\quad\leq~C_\gamma\cdot\bigl(|\mu_1^l|+|\mu_2^l|\bigr)\bigl(]t_0, t_1]\bigr) -C_\gamma\cdot |\mu_1^r|\bigl(]t_0, t_1]\bigr)+\O(1)\cdot\delta_0\,.
\enda
\eeq
 Combining (\ref{Wt01}) with (\ref{we3}) 
 we obtain
\bel{we4}|\mu_1^r|\bigl(]t_0, t_1]\bigr)~\leq~W(t_0)-W(t_1) +3b\, C_\gamma
+\O(1)\cdot\delta_0\,.\eeq
Inserting the bound (\ref{we4}) in (\ref{we1}) we finally obtain
\bel{we5}
|\hat\mu_2| \Big( [\gamma(t_1), y_2[\Big)
~\leq~3b\,C_\gamma + W(t_0) - W(t_1) + 3b\, C_\gamma +\Tilde C\,\delta_0\,,\eeq
where $\Tilde C$ is a constant depending only on $a,b$, 
and $\Tilde K_0$.
We can now choose $\delta_0>0$ small enough so that
\bel{we6}|\hat\mu_2| \Big( [\gamma(t_1), y_2[\Big)
~\leq~7b\,C_\gamma + W(t_0) - W(t_1).
\eeq
Defining
\bel{TK}\Tilde K_{ab} ~\doteq~ 2K_{ab} + 7b C_\gamma\,,\eeq
in view of (\ref{est5}) and (\ref{we6}) we achieve the desired bound (\ref{rar4}).
\endproof
\v
\begin{remark}\label{remark2}{\rm 
The  estimate (\ref{rar1}) remains valid if  the domain $\Gamma$ in 
(\ref{Gd1}) is replaced by
\bel{Gd5}
\Gamma'~=~\Big\{ (t,x)\,;~~t\in [t_0, t_1],\quad
\max\bigl\{  x_1+\hat\lambda\,,~\gamma(t)\bigr\}~\leq~x~\leq~\min
\bigl\{x_2-\hat\lambda
(t-t_0)\,,~\tilde\gamma(t)\bigr\}\Big\},\eeq
where $\gamma(t)$ is the location of a 2-shock, while
$\tilde\gamma(t)$ is the location of a 1-shock (see Fig.~\ref{f:p100}, left).
Indeed,
no wave of any kind can
enter $\Gamma'$ from the left boundary $\gamma$, nor 
from the right boundary $\tilde\gamma$.

Similarly, assume that the domain $\Tilde \Gamma$ in (\ref{Gd3})
is replaced by 
\bel{Gd6}
\Tilde \Gamma'~=~\Big\{ (t,x)\,;~~t\in [t_0, t_1],\quad
\gamma(t)< x< \min \bigl\{  x_2-\hat\lambda
(t-t_0)\,,~\tilde\gamma(t)\bigr\}\Big\},\eeq
where $\tilde\gamma(t) $ is the location of a 1-shock (see Fig.~\ref{f:p100}, right).   Then, setting
$$y_1~\doteq~\gamma(t_1),\qquad\qquad y_2~\doteq~\min  x_2-\hat\lambda
(t_1-t_0)\,,~\tilde\gamma(t_1)\bigr\},$$
the same arguments used in the proof of Lemma~\ref{lemma8}  yield
\bel{rar9} (\hat \mu_1^++\hat \mu_2^+)\bigl([y_1,y_2]\bigr)~\leq~ \Tilde K_{ab}+
\bigl[W(t_0)- W(t_1)\bigr].\eeq
Here $W(t)$ denotes the total weighted strength of waves at time $t$
contained inside $\Tilde\Gamma'$.
Indeed,
no wave of any kind can
enter $\Tilde \Gamma'$ 
from the right boundary $\tilde\gamma$.
}
\end{remark}

\v
\section{Proof of the main theorem}
\label{sec:8}
\setcounter{equation}{0}

\begin{figure}[htbp]
\centering
  \includegraphics[scale=0.5]{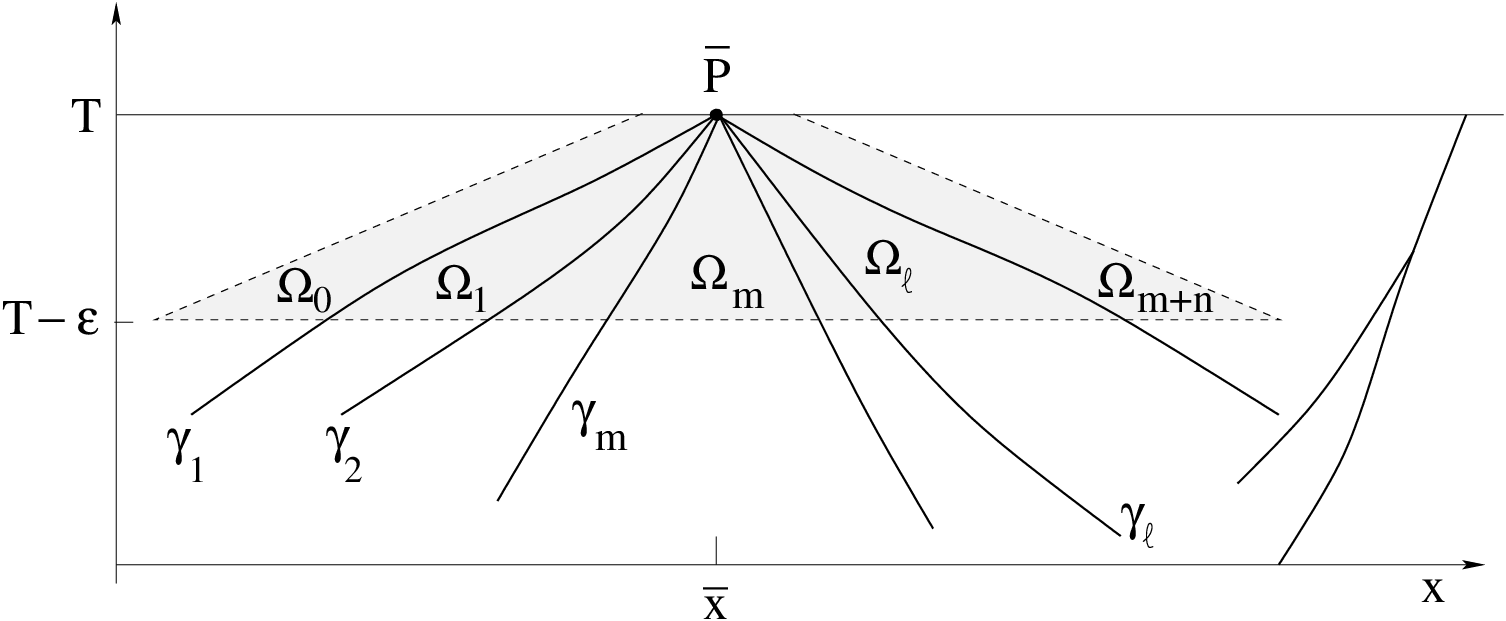}
    \caption{\small By taking a sufficiently small neighborhood
    of $\ov P=(T,\bar x)$ one can assume that all shock curves 
    $\gamma_\ell$ actually meet at $\ov P$.  We illustrate 
    here the case $m=n=3$.  The shaded region is the set $\Gamma_\ve$ in (\ref{Ge}).}
\label{f:p65}
\end{figure}

Toward a proof of Theorem~\ref{thm1}, consider a solution $U=(v,u)$ of (\ref{1}), 
taking values inside the domain $\D_{ab}$, and assume that 
\bel{BU1}
\limsup_{t\to T-} ~\TV\bigl\{ U(t,\cdot)\bigr\}~=~+\infty.\eeq
As remarked in (\ref{away}), blow up of the total variation can occur
only in the region where $|x|\leq R_0 +\hat\lambda T$.
By a compactness argument, there exists a point $\bar x$   such that, for 
every $\ve>0$, 
\bel{BU2}
\limsup_{t\to T-} ~\TV\bigl\{ U(t,\cdot)\,;~~[\bar x-\ve,\, \bar x+\ve]\bigr\}~
=~+\infty.\eeq
If the conclusion of Theorem~\ref{thm1} is violated, then there exists
$\ve>0$ such that the rectangular region 
$$\bigl\{ (t,x)\,;~~t\in [T-\ve,\, T],~~|x-\bar x|\leq\ve\bigr\}$$
contains only finitely many shocks.
By possibly shrinking the value of $\ve$, we can assume that 
all these shocks meet at the point
$\ov P~=~(T,\bar x)$.   In the next steps
of the proof we will derive a contradiction, showing that the 
total variation of $U(t,\cdot)$ in a neighborhood of $\ov P$  
remains uniformly bounded. 
\v
{\bf 1.} 
To fix the ideas, call $x=\gamma_i(t)$, $t\in [t_0,\,T]$ 
the position of the 
$i$-th large shock.  As shown in Fig.~\ref{f:p65} 
we can assume that all these shocks reach 
the point $\ov P$ at time $t=T$, and moreover
\bel{gamn}\gamma_1(t)~<~\cdots~<~\gamma_m(t)~<~\gamma_{m+1}(t)~<~\cdots~<~
\gamma_{m+n}(t),\qquad\qquad T-\ve<t<T.\eeq
Here $\gamma_1,\ldots,\gamma_m$ are 2-shocks while 
$\gamma_{m+1},\ldots,\gamma_{m+n}$ are 1-shocks.
Consider the domain
\bel{Ge}
\Gamma_\ve~\doteq~\Big\{ (t,x)\,;~~t\in [T-\ve,\,T],~~
\bar x-\ve -\hat\lambda(T-t)~<~x~<~\bar x+\ve +\hat\lambda(T-t)\Big\}.\eeq

\begin{figure}[htbp]
\centering
  \includegraphics[scale=0.45]{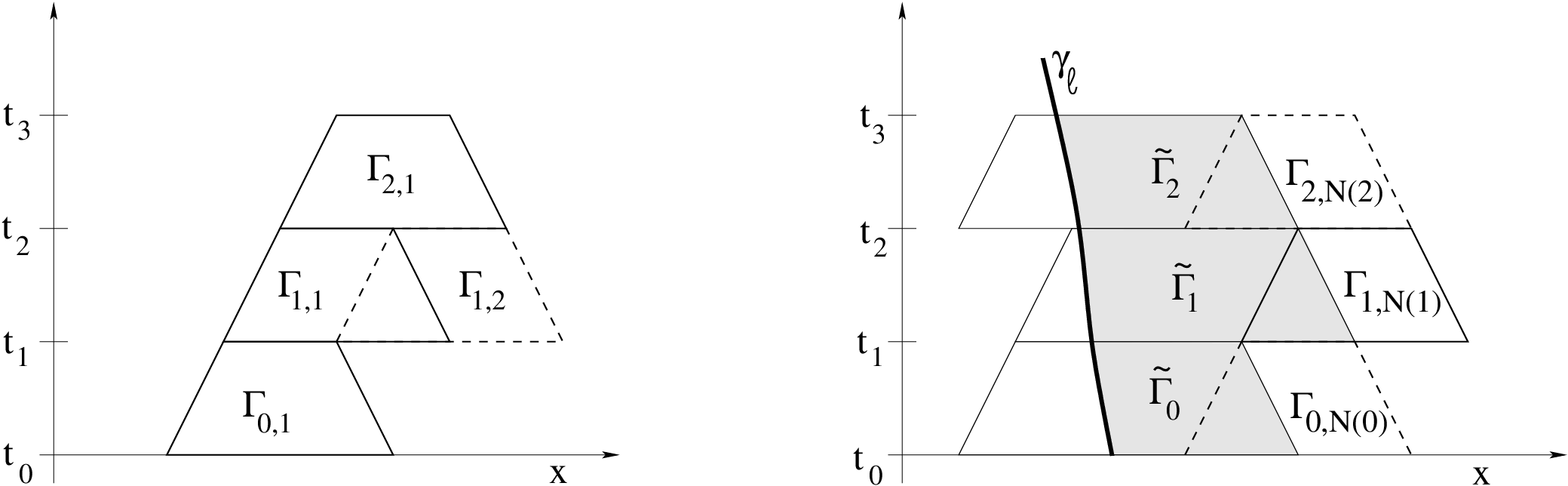}    
  \caption{\small Left: Covering the domain $\Omega_m$ with trapezoids $\Gamma_{j,k}$.  
  Notice that the bottom side of  $\Gamma_{2,1}$
  is the union of the top sides of the two lower trapezoids 
  $\Gamma_{1,1}$ and $\Gamma_{1,2}$. Right:
  Covering the domains $\Omega_\ell$ with trapezoids. 
  Here the B-trapezoids are shaded in grey.
Notice that the bottom side of the boundary trapezoid 
$\Tilde \Gamma_2$ coincides with  
the top side of the $\Tilde \Gamma_1$, while the bottom side of $\Tilde \Gamma_1$
is the union of the top sides of $\Tilde \Gamma_0$ and 
$\Gamma_{0,N(0)}$.
Similarly, the bottom side of $\Gamma_{2,N(2)}$ is 
contained in the union of the top sides
of $\Tilde \Gamma_1$ and $\Gamma_{1,N(1)}$.
}
\label{f:p82}
\end{figure}

Set $t_0\doteq T-\ve$ and 
let $\Omega_0,\Omega_1,\ldots, \Omega_{m+n}$ be the regions 
between these large shocks, so that
\bel{Oj}\Omega_\ell~\doteq~\bigl\{ (t,x)\,;~~\gamma_\ell(t)<x<\gamma_{\ell+1}(t),
\quad t_0<t<T\bigr\}.\eeq
 For convenience, we denote by
$$\gamma_0(t)~\doteq~ \bar x-\ve -\hat\lambda(T-t),
\qquad\qquad \gamma_{m+n+1}(t)~\doteq~\bar x+\ve +\hat\lambda(T-t),$$
 the left and the right side of the trapezoid $\Gamma_\ve$.
By an inductive argument, we will show that the total variation of the solution
restricted to each $\Omega_\ell$ is uniformly bounded.
\v
{\bf 2.} 
To prove that the total variation on the middle domain 
$\Omega_m$ remains bounded, we choose mesh sizes
$\Delta t$, $\Delta x$, with
\bel{Dtx}
\Delta x~=~2\hat\lambda \cdot\Delta t\,,\eeq
As shown in Fig.~\ref{f:p82}, left, we can cover the domain $\Omega_m$ with trapezoids
whose basis has length $\leq 2\Delta x$ and such that the bottom side
of each trapezoid is the union of the top sides of the two lower ones.

More precisely, we set
$$t_j~=~t_0 + j\cdot \Delta t\,,$$ and
cover $\Omega_m$ with trapezoids
of the form (see Fig.~\ref{f:p82}, left)
\bel{Gjk}\bega{l}
\Gamma_{jk}~=~\Omega_m \,\cap\,\Big\{ (t,x)\,;~~t\in [t_j,t_{j+1}]\,,\\[4mm]
\qquad\qquad\qquad k\cdot \Delta x +
j \hat \lambda \cdot\Delta t+ (t-t_j)\hat \lambda<x<(k+2)\cdot \Delta x +
j \hat \lambda \cdot\Delta t - (t-t_j)\hat \lambda\Big\}\,.\enda\eeq
Choosing the mesh $\Delta x>0$ small enough, 
 by Lemma~\ref{lemma3} we can assume that at time $t=t_0$ the total strength
of waves along the bottom side of every trapezoid 
$\Gamma_{0,k}$ is $\leq 2b$.

Let $K_{ab}$ be the constant in Lemma~\ref{lemma4}, 
depending only on $a,b$.
We then choose 
$$K_0 ~=~2(b+2K_{ab})$$ 
and find some 
$\delta_0>0$ so that the conclusion of Lemma~\ref{lemma4} holds.  

We claim that, if all shocks inside $\Omega_m$ have strength 
$\leq\delta_0$, then 
for every $j,k$ one has
 \bel{ws3}
\hbox{[total strength
of waves along the bottom side of $\Gamma_{j,k}\,$]}~\leq~ K_0\,.\eeq
Indeed, since $K_0> 2b$, this is trivially true when $j=0$.

Now assume that (\ref{ws3}) holds for some $j$ and all $k$. 
For a fixed $k$, consider the trapezoid $\Gamma_{j+1,k}$.
Observe that the bottom side of $\Gamma_{j+1,k}$ is the union of the 
top sides of the two lower trapezoids 
$\Gamma_{j,k}$ and $\Gamma_{j, k+1}$. 
By the inductive assumption, an application of Lemma~4 implies that
the total strength of rarefactions contained in the top side of $\Gamma_{j,k}$  
is $\leq K_{ab}$, and the same is true for the 
top side of the trapezoid $\Gamma_{j, k+1}$. 
By (i) in  Lemma~\ref{lemma3}  we conclude that the 
 total strength
of all waves  contained in the bottom side of $\Gamma_{j+1,k}$  
is $\leq 2(b+2K_{ab}) =K_0$.  Hence (\ref{ws3}) holds also 
with $j$ replaced by $j+1$.     By induction, the same estimate
holds for every $j\geq 0$ and every~$k$.

Since $\Omega_m$ can be covered with finitely many trapezoids 
$\Gamma_{j,k}$, this proves that the total variation of the solution
remains uniformly bounded inside $\Omega_m$.
\v
\v
{\bf 3.} Next, assume that the total variation remains bounded 
on the domains $\Omega_m,\ldots, \Omega_{\ell-1}$. We claim that
it also remain bounded on $\Omega_\ell$, the region 
to the right of the 1-shock 
$\gamma_\ell$.

Toward this goal, we consider  the 
trapezoids 
\bel{Gjkl}\bega{l}
\Gamma_{jk}~=~\Omega_\ell \,\cap\,\Big\{ (t,x)\,;~~t\in [t_j,t_{j+1}]\,,\\[4mm]
\qquad\qquad\qquad k\cdot \Delta x +
j \hat \lambda \cdot\Delta t+ (t-t_j)\hat \lambda<x<(k+2)\cdot \Delta x +
j \hat \lambda \cdot\Delta t - (t-t_j)\hat \lambda\Big\}\,.\enda\eeq
but only for values of 
$j,k$ such that $\Gamma_{j,k}$ remains entirely to the right of the 
shock curve 
$\gamma_\ell$, 
namely
$$\gamma_\ell(t_j)~<~k\cdot \Delta x +
j \hat \lambda \cdot\Delta t\,.$$
These will be called {\it I-trapezoids}. 
In addition, for each  $j$ we consider a somewhat wider
trapezoid, of the form
\bel{TGj} \Tilde \Gamma_j~\doteq~
\Big\{ (t,x)\,;~~t\in [t_j,t_{j+1}]\,,~~
\gamma_\ell(t)~<~x~<~(N(j)+1)\cdot \Delta x +j \hat \lambda \cdot\Delta t - (t-t_j)\hat \lambda\Big\}.\eeq
Here  we choose $N(j)$ to be the smallest integer such that
$$\gamma_\ell(t_{j+1})~\leq~N(j)\cdot \Delta x +(j+1) \hat \lambda \cdot\Delta t\,. $$  
%

The $\Tilde \Gamma_j$ will be 
called {\it B-trapezoids}, since they touch the left 
boundary of $\Omega_\ell$.
Since $\gamma_\ell$ is a  1-shock with speed 
$-\hat\lambda<\dot\gamma_\ell<0$, the above choice of 
$N(j)$ guarantees that the lengths of the top and the bottom side of $\Tilde\Gamma_j$ satisfy
\bel{tbG}\left\{\bega{l}
\Delta x~\leq~\bigl[\hbox{length of $\Tilde\Gamma_j^{top}$}\bigr]~\leq~2\Delta x,\\[4mm]
\Delta x~\leq~\bigl[\hbox{length of $\Tilde\Gamma_j^{bottom}$}\bigr]~\leq~3\Delta x.\enda\right.\eeq
 We also observe that
\begi
\item The bottom side of each I-trapezoid $\Gamma_{j,k}$
is contained in the union of the top sides of two
I-trapezoids, or an I-trapezoid and a B-trapezoid.

\item The bottom side of the B-trapezoid $\Tilde \Gamma_j$
is contained in the union of the top side of the lower
B-trapezoid $\Tilde\Gamma_{j-1}$ and (possibly) the top side
of some I-trapezoid.
\endi
By Lemma~\ref{lemma3}, (and the same arguments applied to 
the measures $\mu^l_i$, $i=1,2$  of waves impinging
on the shock $\gamma_\ell$ from the left), 
we can choose the mesh sizes $\Delta t$, $\Delta x$
in (\ref{Dtx}) small enough so that
\begi
\item The total strength of all waves on the bottom side of 
$\Tilde\Gamma_0$ and on the bottom side of each trapezoid 
$\Gamma_{0,k}$ is $\leq 2b$.
\item During every time interval $[t_j,t_{j+1}]$, the total amount of waves impinging on $\gamma_\ell$ from the left
is $\leq 2b$.
\endi
\v
{\bf 4.} Let $K_{ab}$ and $\Tilde K_{ab}$ be the  constants
for which the conclusions  of Lemma~\ref{lemma4} and of Lemma~\ref{lemma8}
hold. Of course, it is not restrictive to assume that $\Tilde K_{ab}> K_{ab}>1$.
Choose 
$$K_0~\doteq~2\bigl(b+ \Tilde K_{ab} + 4(b+K_{ab}) + K_{ab}\bigr)$$
in Lemma~\ref{lemma4} and
$$\Tilde K_0~\doteq~2 K_0 +2(b+2K_{ab})$$
in Lemma~\ref{lemma8}, and let $\delta_0>0$ be a constant 
small enough so that the conclusions of both lemmas hold.

With the above choices, we claim that, for every $j,k$,
\begi
\item[(i)] The total strength of all waves on the bottom side of each 
trapezoid $\Gamma_{jk}$ is $\leq K_0$.
\item[(ii)] The weighted strength of all waves on the bottom side of each 
trapezoid $\Tilde \Gamma_j$ is $W_j(t_j)\leq 2\Tilde K_0$.
\endi

Indeed, by the choice of the step size $\Delta x$,
both claims (i)-(ii) are true at time $t_0=T-\ve$, i.e.~for 
$j=0$ and any $k$.

Arguing by induction, 
assume that the above claims are true for  the trapezoids  $\Gamma_{jk}$
and $\Tilde \Gamma_j$, for a given $j$ and all $k$.   We need to prove that
they hold for the trapezoids   $\Gamma_{j+1, k}$
and $\Tilde \Gamma_{j+1}$ as well.

Call $W_{j}(t)$ the total weighted strength of waves in the trapezoid $\Tilde
\Gamma_{j}$ at time $t\in [t_j, t_{j}]$.
By the inductive assumption $W_j(t_j)\leq 2\Tilde K_{ab}$.
Moreover, by our choice of $\Delta t$ 
the total strength of waves impinging on the 1-shock $\gamma_\ell$
(the left boundary of $\Tilde\Gamma_j$) is $\leq 2b$.
We can thus apply  Lemma~\ref{lemma8} and conclude
\bel{rar44}
[\hbox{total strength of all rarefactions on the upper boundary of 
$\Tilde \Gamma_j$}]~\leq~\Tilde K_{ab} + \bigl[ W_j(t_j)-W_j(t_{j+1})\bigr].\eeq

Two cases must be considered.

CASE 1: $W_j(t_{j})-W_j(t_{j+1})\geq 4(b+K_{ab})$.

Consider the bottom side of $\Tilde \Gamma_{j+1}$.  
As shown in Fig.~\ref{f:p82}  left, this is contained
in the union of the top side of $\Tilde\Gamma_j$ and (possibly)
some smaller trapezoid $\Gamma_{j, N(j)}$.   By Lemma~\ref{lemma4} and 
part (i) of Lemma~\ref{lemma3},   the total strength of waves on the upper
boundary of $\Gamma_{j, N(j)}$ is $\leq 2(b+K_{ab})$.

Hence the total weighted strength of all waves on the lower boundary of
$\Tilde\Gamma_{j+1}$ satisfies the bound
\bel{Wj1}
W_{j+1}(t_{j+1}) ~\leq~W_j(t_{j+1}) + 2 \cdot 2(b+K_{ab})~\leq~W_j(t_j)~\leq~2\Tilde K_0\,.
\eeq

CASE 2: $W_j(t_{j})-W_j(t_{j+1})< 4(b+K_{ab})$.

In this case, by Lemma~\ref{lemma8} the total strength of waves on the upper
boundary of $\Tilde\Gamma_j$ is $< \Tilde K_{ab}+4(b+K_{ab})$.
Since the total strength of positive waves on the upper boundary of the 
small trapezoid $\Gamma_{j, N(j)}$ is $\leq K_{ab}$, 
using again  Lemma~\ref{lemma3} we conclude that the total weighted strength
of all waves on the lower boundary of $\Tilde\Gamma_{j+1}$
satisfies
\bel{Wj2}
W_{j+1}(t_{j+1}) ~\leq~4\Big( b+ \bigl(\Tilde K_{ab}+4(b+K_{ab})\bigr) +
K_{ab}\Big)~\leq~2\Tilde K_0\,.
\eeq
Together, (\ref{Wj1}) and (\ref{Wj2}) yield (ii), with $j$ replaced by $j+1$.

To prove (i), we again consider two cases.

CASE 1: the lower boundary of 
$\Gamma_{j+1,k}$ is contained in the union of the upper boundaries of the two
small trapezoids $\Gamma_{j,k}$ and $\Gamma_{j, k+1}$.

In this case, as in step {\bf 2}, using
Lemma~\ref{lemma4} and Lemma~\ref{lemma3}
we obtain
\bel{Wj3}[\hbox{total strength of waves on the lower boundary of $\Gamma_{j, k+1}$}]
~\leq~2(b + 2K_{ab})~\leq~K_0\,.\eeq
\v
CASE 2: the lower boundary of 
$\Gamma_{j+1,k}$ is contained in the union of the upper boundaries of the 
large trapezoid $\Tilde\Gamma_j$ and of the 
small trapezoid $\Gamma_{j,N(j)}$.

In this case, using
Lemma~\ref{lemma8} and then Lemma~\ref{lemma3}
we obtain
\bel{Wj4}\bega{l}
[\hbox{total strength of waves on the lower boundary of $\Gamma_{j+1, k}$}]
\\[4mm]
\qquad \leq~2(b+W_j(t_{j+1}) + K_{ab})~\leq~2(b+\Tilde K_0 +K_{ab})~\leq~
K_0\,.\enda\eeq
Together, (\ref{Wj3}) and (\ref{Wj4}) yield (i), with $j$ replaced by $j+1$.

By induction on $j$, both of our claims are thus proved.
Since the domain $\Omega_\ell$ can be covered with finitely many 
trapezoids $\Tilde\Gamma_j$ or $\Gamma_{jk}$, we conclude that the
total variation of the solution $U=U(t,x)$ remains uniformly bounded, 
restricted to the domain $\Omega_\ell$.
\v
{\bf 5.}
The proof of Theorem~\ref{thm1} is 
now achieved by induction on $\ell=m+1,\ldots, m+n$.
The analysis of the total variation on the domains
$\Omega_{m-1},\ldots, \Omega_0$ is entirely similar.
\endproof
\v

{\bf Acknowledgment.} The research of the first author
was partially supported by NSF, with grant DMS-1411786: Hyperbolic Conservation Laws and Applications.
The research of the second author was partially 
supported by NSF with grant DMS-1715012.

\v


\begin{thebibliography}{99}
\bibitem{AS}  D.~Amadori and W.~Shen, Global existence of large BV solutions in a model of granular flow. 
{\it Comm. Partial Differential Equations}
{\bf  34} (2009), 1003--1040.

\bibitem{BJ}
P.~Baiti and H.~K.~Jenssen,
Blowup in $\L^\infty$ for a class of genuinely nonlinear hyperbolic systems of conservation laws.
{\it Discrete Contin. Dynam. Systems} {\bf 7} (2001), 837--853.

\bibitem{BCM} S.~Bianchini, R.~M.~Colombo, and F.~Monti,
$L^\infty$  solutions for $2\times 2$ systems of conservation laws.
{\it Riv. Math. Univ. Parma} {\bf 1} (2010), 189--204. 

\bibitem{Bbook}
A.~Bressan, {\it Hyperbolic Systems of Conservation Laws. The One
Dimensional Cauchy Problem}. Oxford University Press, 2000.

\bibitem{BCZ} A.~Bressan, G.~Chen, Q.~Zhang, 
Lack of BV bounds for approximate solutions
to the p-system with large data, {\it J. Differential Equations} 
{\bf 256} (2014), 3067--3085.

\bibitem{BCZZ} A.~Bressan, G.~Chen, Q.~Zhang, and S.~Zhu,
 No BV bounds for approximate solutions to p-system with 
general pressure law, {\it J. Hyperbolic Diff. Equat.} {\bf 12} (2015), 1--18.

\bibitem{BC1}  A.~Bressan and R.~M.~Colombo,
Unique solutions of $2\times 2$ conservation laws with large
data, {\it Indiana Univ. Math. J.} {\bf 44} (1995),
677--725.


\bibitem{BC2}  A.~Bressan and R.~M.~Colombo,
Decay of positive waves in nonlinear systems of conservation
laws, {\it Ann. Scuola Normale Superiore Pisa} {\bf IV - 26} (1998),
133--160.

\bibitem{BLF}
A.~Bressan and P.~LeFloch, Structural stability and regularity 
of entropy solutions to hyperbolic systems of conservation laws,
{\it Indiana Univ. Math. J.} {\bf 48} (1999), 43--84.

\bibitem{BLY} A.~Bressan, T.~P.~Liu and T.~Yang,  $ L^1$ stability
estimates for $n\times n$ conservation laws,
{\it Arch. Rational Mech. Anal.} {\bf 149} (1999), 1--22.

\bibitem{BY} A.~Bressan and  T.~Yang,
A sharp decay estimate for positive nonlinear waves,
 {\it SIAM Jour. Math. Anal.} {\bf 36} (2004), 659-677.
 

\bibitem{CH} T.~Chang and L.~Hsiao,
   {\it The Riemann problem and interaction of waves in gas dynamics},
   Longman Scientific \& Technical, Harlow, 1989.
   
\bibitem{CG} G.~Chen, {Optimal time-dependent lower bound on density for
classical solutions of 1-D compressible Euler equations}, 
{\it Indiana Univ. Math. J. } {\bf 66} (2017), 725--740.

\bibitem{Daf}
C.~Dafermos, Generalized characteristics and the structure of 
solutions of hyperbolic conservation laws, {\it Indiana Univ. Math. J.} {\bf 26} (1977), 1097--1119.

\bibitem{DP} R.~DiPerna, Singularities of solutions of nonlinear hyperbolic 
systems of conservation laws, {\it Arch. Rational Mech. Anal.} {\bf 60}
(1975), 75--100.


\bibitem{CJ}
G.~Chen and H.~K.~Jenssen, No TVD fields for 1-d isentropic gas flow,
{\it Comm. Partial Differential Equations}, {\bf 38} (2013), 629--657.

\bibitem{C}  C.~Cheverry,  Syst\`emes de lois de conservation et 
stabilit\'e BV.  {\it Mem. Soc. Math. France} {\bf 75} (1998).

\bibitem{CF} R.~Courant and K.~O.~Friedrichs, {\it Supersonic flow
and shock waves}, Wiley-Interscience, New York, 1948.

\bibitem{G} J.~Glimm, Solutions in the large for nonlinear hyperbolic systems
of equations, {\it Comm. Pure Appl. Math.} {\bf 18} (1965), 697--715.

\bibitem{GL} J.~Glimm and P.~Lax, Decay of solutions of systems of
nonlinear hyperbolic
conservation laws, {\it Amer. Math. Soc. Memoir} {\bf 101} (1970).


\bibitem{Hoff} D.~Hoff,
Invariant regions for systems of conservation laws.
{\it Trans. Amer. Math. Soc.} {\bf  289} (1985), 591--610.

\bibitem{HR}
H.~Holden and N.~H.~Risebro,
{\it Front Tracking for Hyperbolic Conservation Laws.}  Springer-Verlag,
 New York, 2002.

\bibitem{J} H.~K.~Jenssen, Blowup for systems of conservation laws,
{\it SIAM J. Math. Anal.} {\bf 31} (2000), 894--908.


\bibitem{Lw}
M.~Lewicka, Well-posedness for hyperbolic systems of conservation laws with large BV data.
{\it Arch. Rational Mech. Anal.}
{\bf   173}  (2004),  415--445.


\bibitem{Liu}
T.~P.~Liu,  Linear and nonlinear large-time behavior of solutions of general systems of hyperbolic conservation laws. {\it Comm. Pure Appl. Math.} {\bf 30} (1977), 767--796. 

\bibitem{MS} C.~Moler and J.~Smoller, Elementary interactions in quasi-linear hyperbolic systems, {\it Arch.
Rational Mech. Anal.} {\bf 37} (1970), 309--322.

\bibitem{N}
T.~Nishida,
Global solution for an initial boundary value problem of a quasilinear hyperbolic system.
{\it Proc. Japan Acad.} {\bf 44} (1968), 642--646.

\bibitem{NS}
T.~Nishida and J.~Smoller,
Solutions in the large for some nonlinear hyperbolic conservation laws. 
{\it Comm. Pure Appl. Math.} {\bf 26} (1973), 183--200. 


\bibitem{Sm} J.~Smoller, {\it Shock waves and reaction-diffusion equations}, Second edition.
  Springer-Verlag, New York, 1994.

\bibitem{TY} B.~Temple and R.~Young. The large time stability of sound waves. {\it Comm. Math. Phys.} {\bf 179} (1996), 417--466. 

\bibitem{Y} R.~Young, Sup-norm stability for Glimm's scheme. 
{\it Comm. Pure Appl. Math.} {\bf  46} (1993), 903--948. 

\end{thebibliography}
\end{document}